\documentclass[12p;t]{article}
\usepackage{ucs}
\usepackage[utf8x]{inputenc} 
\usepackage[russian]{babel}  
\usepackage{amsfonts}
\usepackage{amssymb}
\usepackage{amsmath}
\usepackage{euscript}
\newcommand{\EHR}{\mathrm{EHR}}

\title{Spectrum for first-order properties of random hypergraphs} 
\author{}
\date{}
\begin{document}
\maketitle

\begin{center}
{\large S.N. Popova~\footnote{Moscow Institute of Physics and Technology, Laboratory of Advanced Combinatorics and Network Applications.}}
\end{center}

The notion of spectrum for first-order properties introduced by J. Spencer for Erdős–Rényi random graph is considered in relation to random uniform hypergraphs.
We study properties of spectrum for first-order formulae with bounded quantifier depth and estimate the values of minimum and maximum points of this spectrum. 
We also consider the set of limit points of the spectrum and give a bound for its minimum value.


\section{Introduction}
\indent 

Asymptotic behavior of first-order properties probabilities of the Erdős–Rényi random graph
$G(n, p)$ has been widely studied (see \cite{shelah_spencer}--\cite{raigor}).
In \cite{inf_spectra} the notion of spectrum of first-order properties was introduced and it was shown that there exists a first-order property with infinite spectrum.
In this paper, we consider spectrum of first-order properties in relation to random uniform hypergraphs.
Let us define the \textit{random $s$-uniform hypergraph} $G^s(n, p)$. Consider the set $\Omega_n = \{G = (V_n, E)\}$ of all $s$-uniform hypergraphs ($s$-hypergraphs)
with the set of vertices $V_n = \{1, 2, \ldots, n\}$. The random hypergraph $G^s(n, p)$ is a random element with probability distribution
$$\Pr[G^s(n, p) = G] = p^{|E|} (1 - p)^{{n \choose s} - |E|}.$$
Note that for $s = 2$ in this definition we obtain the Erdős–Rényi random graph $G(n, p)$.
Let us denote the event ``$G^s(n, p)$ has a property $L$'' by $G^s(n, p) \models L$.

\textit{First-order properties} of $s$-uniform hypergraphs are defined by first-order formulae (see \cite{raigor}, \cite{ebbing}) which are built of predicate symbols $N, =$, logical connectivities, variables and quantifiers $\forall, \exists$. The $s$-ary predicate symbol $N$ expresses the property of $s$ vertices to constitute an edge.  
Let us remind that the \textit{quantifier depth} (see \cite{raigor}, \cite{ebbing}) of a first-order formula is the maximum number of nested quantifiers.
Let $\mathcal L_k$ denote the set of all first-order properties which can be expressed by first-order formulae with quantifier depth at most $k$.

We say that the random hypergraph $G^s(n, p)$ \textit{obeys the zero-one law} if for any first-order property $L$ the probability $\Pr[G^s(n, p) \models L]$ tends either to 0 or to 1 as $n \to \infty$. We say that the random hypergraph $G^s(n, p)$ \textit{obeys the zero-one $k$-law} if for any first-order property $L \in \mathcal L_k$ the probability $\Pr[G^s(n, p) \models L]$ tends either to 0 or to 1 as $n \to \infty$.

For any first order property $L$ we define its \textit{spectrum} $S(L)$. $S(L)$ is the set of all $ \alpha \in (0, s - 1)$ 
(we take the interval $(0, s - 1)$ because the structure of $G^s(n, n^{-\alpha})$ for $\alpha > s - 1$  is considerably simpler than that for $\alpha < s- 1$ and the study of $G^s(n, p)$ with $p(n) \gg n^{- s + 1}$ is a subject of a separate investigation) which do not satisfy the
following property: $\lim_{n \to \infty} \Pr[G^s(n, p) \models L]$ exists and is either $0$ or $1$.
Denote the union of $S(L)$ over all $L\in \mathcal L_k$ by $S_k$.
In other words, $S_k$ is the set of all $\alpha \in (0, s - 1)$ for which the random hypergraph $G^s(n, n^{-\alpha})$ does not obey the zero-one $k$-law. 
Let $(S_k)'$ be the set of limit points in $S_k$.

S. Shelah and J. Spencer (see \cite{shelah_spencer}) showed that when $\alpha$ is an irrational number and $p(n) = n^{-\alpha + o(1)}$ then the random graph $G(n, p)$ obeys the zero-one law. If $\alpha \in (0, 1) \cap \mathbb Q$, then $G(n, n^{-\alpha})$ does not obey the zero-one law. In \cite{matushkin} an extension of the result from \cite{shelah_spencer} to the case of random uniform hypergraph $G^s(n, p)$ was established.

M. E. Zhukovskii proved the zero-one $k$-law for the random graph $G(n, n^{-\alpha})$ where $\alpha$ is close to 0 or close to 1. He also found the minimum and maximum points in $S_k$.

\textbf{Theorem 1} (\cite{zhuk1}).
\textit{
If $\alpha < \frac{1}{k-2}$, then the random graph $G(n,n^{-\alpha})$ obeys the zero-one $k$-law. If $\alpha = \frac{1}{k-2}$, then $G(n,n^{-\alpha})$ does not obey the zero-one $k$-law.
}\\

\textbf{Theorem 2} (\cite{zhuk2},~\cite{zhuk3}).
\textit{
Let $k>3$ be an arbitrary natural number. Let $\mathcal{Q}$ be the set of all positive rational numbers with numerator at most $2^{k-1}$, and let $\tilde{\mathcal{Q}}$ be the set of natural numbers not exceeding $2^{k-1}-2$. The random graph $G(n,n^{-\alpha})$ obeys the zero-one $k$-law if $\alpha = 1 - \frac{1}{2^{k-1}+\beta}$, where $\beta \in (0,\infty) \setminus \mathcal{Q}$. If $\alpha = 1 - \frac{1}{2^{k-1}+\beta}$, $\beta \in \tilde{\mathcal{Q}}$, then $G(n,n^{-\alpha})$ does not obey the zero-one $k$-law. 
If $\alpha \in  \{1 - \frac{1}{2^k}, 1 - \frac{1}{2^k - 1}\}$, then $G(n,n^{-\alpha})$ obeys the zero-one $k$-law. 
}\\

Theorems 1 and 2 imply that $\min S_k = \frac{1}{k - 2}$ and $\max S_k = 1 - \frac{1}{2^k  - 2}$.

M. E. Zhukovskii and J. Spencer examined the question about the minimum and the maximum limit points in $S_k$ for the random graph $G(n, n^{-\alpha})$.

\textbf{Theorem 3} (\cite{zhuk_spencer}).
\textit{
There exists such $k_0$ such that for any natural $k > k_0$ 
$$\min (S_k)' \le \frac{1}{k - 11}, \quad \max (S_k)' \ge 1 - \frac{1}{2^{k - 13}}.$$
}

M. E. Zhukovskii obtained an upper bound for the minimum limit point in $S_k$ which is better than the bound from theorem 3 for small $k$.

\textbf{Theorem 4} (\cite{zhuk4}).
\textit{
For any $k \ge 5$, $\frac{1}{\lfloor k/2 \rfloor} \in (S_k)'$.
}\\

Theorem 4 implies the following corollary about the minimum quantifier depth $k$ for which $S_k$ is infinite.

\textbf{Corollary 1} (\cite{zhuk_form}).
\textit{
The minimum $k$ such that the set $S_k$ is infinite equals 4 or 5.
}\\

Note that since $|\mathcal L_k| < \infty$ (see Proposition 3.1.3 in \cite{ebbing}), the set $S_k$ is infinite if and only if there exists $L \in \mathcal L_k$ with infinite spectrum $S(L)$.
The following theorem gives an almost exact bound for the minimum quantifier depth of a first-order formula with an infinite spectrum in the case of $s$-uniform hypergraphs, where $s \ge 3$.

\textbf{Theorem 5} (\cite{infspectrum}). \textit{For any $s \ge 3$, the minimum $k$ such that the set $S_k$ is infinite equals $s + 1$ or $s + 2$.}\\

\section{New results}
\indent

We found an interval with the left endpoint equals to 0, such that for all $\alpha$ from this interval the random hypergraph $G^{s}(n,n^{-\alpha})$ obeys the zero-one $k$-law. Note that 
$G^{s}(n,n^{-\alpha})$ obeys the zero-one $s$-law for all $\alpha \in (0, s - 1)$ so we study zero-one $k$-laws for $k \ge s + 1$. 

\textbf{Theorem 6.} \textit{
Let $s \ge 3$, $k \ge s + 1$, $\alpha > 0$ and 
$$\frac{1}{\alpha} > {{k - 1} \choose {s - 1}} - 1 - \frac{s - 1}{k - 1} + \frac{2\left(1 + \frac{s - 1}{k - 1}\right)}{{{k - 1} \choose {s - 1}} + 2}.$$
Then the random hypergraph $G^s(n, n^{-\alpha})$ obeys the zero-one $k$-law.
}\\

We also prove that near the right endpoint of our interval there is a point $\alpha$ for which $G^{s}(n,n^{-\alpha})$ does not obey the zero-one $k$-law.

\textbf{Theorem 7.} \textit{
Let $s \ge 3$, $k \ge s + 2$. Then there exists $\alpha > 0$ such that 
$$\frac{1}{\alpha} > {{k - 1} \choose {s - 1}} - 1 - \frac{s - 1}{k - 1} - \frac{2}{{{k - 1} \choose {s - 1}}}$$
and $G^s(n, n^{-\alpha})$ does not obey the zero-one $k$-law.
}\\

The question what is the exact value of the minimal point for which the zero-one $k$-law does not hold remains open. Theorems 6 and 7 imply that this value is very close to
$1/\left({{k - 1} \choose {s - 1}} - 1 - \frac{s - 1}{k - 1}\right)$. 

We also examine the zero-one $k$-law for the left neighborhood of $s-1$. Our results are the following. 

\textbf{Theorem 8}. \textit{
Let $\alpha \in \left(s - 1 - \frac{1}{2^{k - s + 1}}, s - 1\right) \setminus \mathcal Q_k$, where
$\mathcal Q_k = \{s - 1 - \frac{1}{2^{k - s + 1} + \frac{a}{b}} | a, b \in \mathbb N, a \le 2^{k - s + 1}\}$.
Then the random hypergraph $G^s(n, n^{-\alpha})$ obeys the zero-one $k$-law. 
}\\

\textbf{Theorem 9}. \textit{
Let $s \ge 3$, $k \ge s + 4$, $\alpha = s - 1 - \frac{1}{2^{k - s + 1} + a}$, where $a \in \mathbb N$, $a \le 2^{k - s - 2} + 2^{k - s - 3} + 1$.
Then the random hypergraph $G^s(n, n^{-\alpha})$ does not obey the zero-one $k$-law.
}\\

For an $s$-hypergraph $G$, denote by $d_G(x, y)$ the distance between vertices 
$x, y \in V(G)$ which is defined as the minimal length of a path between $x$ and $y$ (the length of a path is the number of its edges).

Theorem 9 gives a weaker result than theorem 2 which is related to the fact that
 the following property (which was used for the proof of theorem 2) is true for graphs but is not true for $s$-hypergraphs, where $s \ge 3$: if for a graph $G$ there exist vertices $x_1, x_2, x_3 \in V(G)$ such that $d_G(x_1, x_2) = a$, $d_G(x_1, x_3) = \lfloor (a + 1)/2 \rfloor$, $d_G(x_2, x_3) = \lceil (a + 1)/2 \rceil$ then there exists a cyclic $(2a + 1)$-extension in $G$.

We also study the question about the minimum limit point in $S_k$.

\textbf{Theorem 10.} \textit{
Let $k \ge s + C$, where $C$ is a constant. Then $\min (S_k)' \le \frac{1}{{{k - 11} \choose {s - 1}}}$.
}\\

Taking into account the bound from theorem 6 we deduce that theorem 10 gives an asymptotically tight bound for $\min (S_k)'$.

\textbf{Corollary 2.} \textit{
We have $\min (S_k)' \sim \frac{(s-1)!}{k^{s - 1}}$, as $k \to \infty$.
}\\

We also obtain a bound for $\min (S_k)'$ which is better than the bound from theorem 10 for small $k$.

\textbf{Theorem 11.} \textit{
Let either $s = 2$ and $k \ge 5$ or $s \ge 3$ and $k \ge s + 2$. Then $\frac{1}{{l(k) \choose {s - 1}}} \in (S_k)'$, where $l(k) = \max\{l \in \mathbb N: {l \choose {s - 1}} (l + 2) \le {k \choose s}\}$.
}\\

This theorem gives a bound for the minimum limit point in $S_k$ for all $k$ for which it is known that $S_k$ is infinite (see corollary 1 and theorem 5).
Note that if $s = 2$, then $l(k) = \left\lfloor \sqrt{k(k - 1)/2 - 1} \right\rfloor - 1$ and the bound for $\min (S_k)'$ obtained in theorem 11 improves the bound from theorem 4.

\section{Auxiliary statements}

\subsection{Small sub-hypergraphs and extensions}
\indent

For an arbitrary $s$-hypergraph $G = (V, E)$, set $v(G) = |V|$, $e(G) = |E|$, $\rho(G) = \frac{e(G)}{v(G)}$, $\rho^{\max}(G) = \max_{H \subseteq G} \rho(H)$. $G$ is called \textit{strictly balanced} if the \textit{density} $\rho(G)$ of this graph is greater than the density of any its proper subhypergraph.
Denote the property of containing a copy of $G$ by $L_G$.

\textbf{Theorem 12} (\cite{vantsyan}). \textit{If $p \ll n^{-1/\rho^{\max}(G)}$, then
$\lim_{n \to \infty} \Pr[G^s(n, p) \models L_G] = 0$.
If $p \gg n^{-1/\rho^{\max}(G)}$, then
$\lim_{n \to \infty} \Pr[G^s(n, p) \models L_G] = 1$.}\\

In other words, the function $n^{-1/\rho^{\max}(G)}$ is a threshold for the property $L_G$.

Let $G_1, \dots, G_m$ be strictly balanced $s$-hypergraphs, $\rho(G_1) = \ldots = \rho(G_m) = \rho$. Let $a_i$ be the number of automorphisms of $G_i$.
Denote by $N_{G_i}$ the number of copies of $G_i$ in $G^{s}(n,p)$. The following theorem is a generalization of a classical result of Bollob\'{a}s (see, for example,~\cite{JLR}) to the case of $s$-uniform hypergraphs.

\textbf{Theorem 13}. \textit{
If $p=n^{-1/\rho}$, then
$$
(N_{G_1}, \dots, N_{G_m}) \xrightarrow{d}(P_1, \dots, P_m),
$$
where $P_i \sim \mathrm{Pois}\left(\frac{1}{a_i}\right)$ are independent Poisson random
variables.
}\\

Consider arbitrary $s$-hypergraphs $G$ and $H$ such that $H \subset G$, $V(H) = \{x_1, \ldots, x_l\}$, $V(G) = \{x_1, \ldots, x_m\}$. Denote $v(G, H) = v(G) - v(H)$, $e(G, H) = e(G) - e(H)$,  $\rho(G, H) = \frac{e(G, H)}{v(G, H)}$,
$\rho^{\max}(G, H) = \max_{H \subset K \subseteq G} \rho(K, H)$. 
The pair $(G, H)$ is called \textit{strictly balanced} if $\rho(G, H) > \rho(K, H)$ for any graph $K$ such that
$H \subset K \subset G$.

Consider $s$-hypergraphs $\tilde H, \tilde G$, where $V(\tilde H) = \{\tilde x_1, \ldots, \tilde x_l\}$, $V(\tilde G) = \{\tilde x_1, \ldots, \tilde x_m\}$,
$\tilde H \subset \tilde G$. The hypergraph $\tilde G$ is called \textit{$(G, (x_1, \ldots, x_l))$-extension} of the collection $(\tilde x_1, \ldots, \tilde x_l)$, if
$$
\{x_{i_1}, \ldots, x_{i_s}\} \in E(G) \setminus E(H) \Rightarrow \{\tilde x_{i_1}, \ldots, \tilde x_{i_s}\} \in E(\tilde G) \setminus E(\tilde H).
$$ 
If
$$
\{x_{i_1}, \ldots, x_{i_s}\} \in E(G) \setminus E(H) \Leftrightarrow \{\tilde x_{i_1}, \ldots, \tilde x_{i_s}\} \in E(\tilde G) \setminus E(\tilde H),
$$ 
we call $\tilde G$ a \textit{strict $(G, (x_1, \ldots, x_l))$-extension} of $(\tilde x_1, \ldots, \tilde x_l)$.

Let $\alpha > 0$. For any pair $(G, H)$, where $H \subset G$, set $f_\alpha(G, H) = v(G, H) - \alpha e(G, H)$. 
The pair $(G, H)$ is called \textit{$\alpha$-safe}, if
$f_\alpha(K, H) > 0$ for any $K$, $H \subset K \subseteq G$. If for any $K$ such that $H \subset K \subseteq G$ we have
$f_\alpha(G,K) < 0$, then the pair $(G, H)$ is called \textit{$\alpha$-rigid}. If for any $K$ such that $H \subset K \subset G$ we have
$f_\alpha(K, H) > 0$ and $f_\alpha(G, H) = 0$, then the pair $(G, H)$ is called \textit{$\alpha$-neutral}.

Let $\tilde H \subset \tilde G \subset \Gamma$, $T \subset K$ and $|V(T)| \le |V(\tilde G)|$. The pair $(\tilde G, \tilde H)$ is called \textit{$(K, T)$-maximal} in $\Gamma$, if for any subhypergraph $\tilde T \subset \tilde G$ with $|V(\tilde T)| = |V(T)|$ and $\tilde T \cap \tilde H \neq \tilde T$ there does not exist
a strict $(K, T)$-extension $\tilde K$ of $\tilde T$ in the hypergraph $\Gamma \setminus (\tilde G \setminus \tilde T)$ such that $E((\tilde K \cup \tilde G) \setminus \tilde T) \setminus (E(\tilde K \setminus \tilde T) \cup
E(\tilde G \setminus \tilde T)) = \varnothing$.
The hypergraph $\tilde G$ is called \textit{$(K, T)$-maximal} in $\Gamma$, if for any subhypergraph $\tilde T \subset \tilde G$ with $|V(\tilde T)| = |V(T)|$
there does not exist
a strict $(K, T)$-extension $\tilde K$ of $\tilde T$ in the hypergraph $\Gamma \setminus (\tilde G \setminus \tilde T)$ such that $E((\tilde K \cup \tilde G) \setminus \tilde T) \setminus (E(\tilde K \setminus \tilde T) \cup
E(\tilde G \setminus \tilde T)) = \varnothing$.

Let a pair $(G, H)$ be $\alpha$-safe and $\mathcal K_r$ be the set of all $\alpha$-rigid and $\alpha$-neutral pairs $(K, T)$, where
$|V(T)| \le |V(G)|$, $|V(K) \setminus V(T)| \le r$. Let $\tilde x_1, \ldots, \tilde x_l \in V_n$.
Define a random variable $N^{r}_{(G, H)}(\tilde x_1, \ldots, \tilde x_l)$ which assigns to each hypergraph $\mathcal G \in \Omega_n$ the number of all strict
$(G, H)$-extensions $\tilde G$ of the hypergraph $\tilde H = \mathcal G|_{\{\tilde x_1, \ldots, \tilde x_l\}}$ such that for any $(K, T) \in \mathcal K_r$ the pair $(\tilde G, \tilde H)$ is $(K, T)$-maximal in $\mathcal G$.

The following theorem is the generalization of Lemma 10.7.6 from \cite{alon_spencer} and Proposition 1 from \cite{zhuk_ext} to the case of $s$-uniform hypergraphs.

\textbf{Theorem  14.}
\textit{
Let $p(n) = n^{-\alpha}$, $r \in \mathbb N$ and $(G, H)$ be $\alpha$-safe.
Then a.a.s. for every vertices $\tilde x_1, \ldots, \tilde x_l$ the following relation holds:
$$
N^{r}_{(G, H)}(\tilde x_1, \ldots, \tilde x_l) \sim {\sf E}[N^{r}_{(G, H)}(\tilde x_1, \ldots, \tilde x_l)] = \Theta(n^{f_\alpha(G, H)}).
$$
}\\

Denote by $L^{r}_{(G, H)}$ the following property: for any collection $(\tilde x_1, \ldots, \tilde x_l)$ there exists a strict $(G, (x_1, \ldots, x_l))$-extension which is $(K, T)$-maximal for all $(K, T) \in \mathcal K_r$. 
Theorem 14 implies that if $(G, H)$ is $\alpha$-safe, then the random hypergraph $G^s(n, n^{-\alpha})$ satisfies $L^r_{(G, H)}$ a.a.s.

Denote by $\tilde L_{(G, H)}$ the property that there exists a copy of $H$ such that no copy of $G$ contains it.
Denote by $\tilde N_{(G, H)}$ the number of copies of $H$ in $G^s(n, p)$ which are not contained in any copy of $G$.

\textbf{Proposition 1} (\cite{infspectrum}). \textit{Let $H$ be a strictly balanced $s$-hypergraph, $(G, H)$ be a strictly balanced pair, $\rho(H) = \rho(G, H) = 1/\alpha$.
Let $a_1$ be the number of automorphisms of $H$ which are extendable to some automorphism of $G$. Let $a_2$ be the number of automorphisms 
$\sigma: V(G) \to V(G)$ with $\sigma(x) = x$ for all $x \in V(H)$.
Then $$\tilde N_{(G, H)} \stackrel{d}{\longrightarrow} \mathrm{Pois}\left(\frac{1}{a(H)}\exp\left(-\frac{a(H)}{a_1 a_2}\right)\right).$$}\\

\subsection{Cyclic extensions}
\indent

Consider a pair of $s$-hypergraphs $(G, H)$, $G \supset H$. We say that $G$ is a \textit{cyclic $m$-extension} of $H$ if one of the following conditions holds:
\begin{itemize}
\item there exists a vertex $x_1 \in V(H)$ such that
$$V(G) \setminus V(H) = \{y_1, \ldots, y_{k(s - 1)}, z_1, \ldots, z_l\},$$
\begin{multline*}
E(G) \setminus E(H) = \{\{x_1, y_1, \ldots, y_{s - 1}\}, \ldots, \{y_{(k - 1)(s - 1)}, \ldots, y_{k(s - 1)}\}, \\
\{y_{k(s - 1)}, z_1, \ldots, z_l, u_1, \ldots, u_{s - 1 - l}\}\},
\end{multline*}  
where $1 \le k \le m - 1$, $0 \le l < s - 1$ and $u_1, \ldots, u_{s - 1 - l} \in \{y_1, \ldots, y_{k(s - 1) - 1}\}$ are distinct,
\item there exist distinct vertices $x_1, x_2 \in V(H)$ such that
$$V(G) \setminus V(H) = \{y_1, \ldots, y_{k(s - 1)}, z_1, \ldots, z_l\},$$
\begin{multline*}
E(G) \setminus E(H) = \{\{x_1, y_1, \ldots, y_{s - 1}\}, \ldots, \{y_{(k - 1)(s - 1)}, \ldots, y_{k(s - 1)}\}, \\
\{x_2, y_{k(s - 1)}, z_1, \ldots, z_l, u_1, \ldots, u_{s - 2 - l}\}\},
\end{multline*}
where $1 \le k \le m - 1$, $0 \le l \le s - 2$ and $u_1, \ldots, u_{s - 2 - l} \in \{y_1, \ldots, y_{k(s - 1) - 1}\}$ are distinct,
\item  there exist distinct vertices $x_1, \ldots, x_l \in V(H)$, where $2 \le l \le s - 1$, such that
$$V(G) \setminus V(H) = \{y_1, \ldots, y_{s - l}\}, E(G) \setminus E(H) = \{x_1, \ldots, x_l, y_1, \ldots, y_{s - l}\},$$
\end{itemize}
and  $\rho^{\max}(G) < \frac{m}{m(s - 1) - 1}$.

For any $m \ge 1$, let us define a set of $s$-hypergraphs $\mathcal H_m$. The hypergraph $(\{x\}, \varnothing)$, where $x$ is a vertex, belongs to $\mathcal H_m$. If $H \in \mathcal H_m$, then $\mathcal H_m$ contains all cyclic $m$-extensions of $H$ and all $s$-hypergraphs $G$ such that $V(G) = V(H)$, $E(G) \supset E(H)$ and $\rho^{\max}(G) < \frac{m}{m(s - 1) - 1}$.

Note that for any hypergraph $G \in \mathcal H_m$ there exists a sequence of hypergraphs 
$G_0 = (\{x\}, \varnothing) \subsetneq G_1 \ldots \subsetneq G_t \subseteq G$ such that $G_{i + 1}$ is a cyclic $m$-extension of $G_i$ for all $i \in \{0, \ldots, t - 1\}$. Let us call such sequence of hypergraphs \textit{$m$-decomposition} of $G$. 

\textbf{Proposition 2}.
\textit{
\begin{itemize}
\item[1)] Let $G \in \mathcal H_m$, $G \neq (\{x\}, \varnothing)$. Then either $G$ is a cyclic $m$-extension of $(\{x\}, \varnothing)$ and $\frac{1}{\rho^{\max}(G)} = s - 1$ or there exist 
$a, b \in \mathbb N$ such that $a \le m$ and $\frac{1}{\rho^{\max}(G)} = s - 1 - \frac{1}{m + a/b}$.
\item[2)] Let $m \in \mathbb N$, $\rho < \frac{m}{m(s - 1) - 1}$. Then there exists $\eta = \eta(\rho) \in \mathbb N$ such that
any hypergraph $G \in \mathcal H_m$ with $|V(G)| \ge \eta$ contains a subhypergraph $H \in \mathcal H_m$ with $|V(H)| \le \eta$ and density
$\rho(H) > \rho$.
\end{itemize}
}

Let $n_1, n_2, n_3, n_4 \in \mathbb N$, $n_2 \le n_1$, $n_4 \le n_3$, $\rho > 0$.  Denote by $\mathcal K^{\rho}$ the set of all $1/\rho$-rigid pairs $(K_1, K_2)$ with $v(K_1) \le n_3$, $v(K_2) \le n_4$.We say that a $s$-hypergraph $G$ is
$(n_1, n_2, n_3, n_4, \rho)-sparse$ if it posseses the following properties.
\begin{itemize}
\item[1)] Let $H$ be a $s$-hypergraph with at most $n_1$ vertices. If $\rho^{\max}(H) > \rho$, then $G$ does not contain a subhypergraph isomorphic to $H$. 
If $\rho^{\max}(H) < \rho$, then $G$ contains an induced subhypergraph which is isomorphic to $H$ and $(K_1, K_2)$-maximal in $G$ for all $(K_1, K_2) \in \mathcal K^{\rho}$.  
\item[2)] For any $1/\rho$-safe pair $(H_1, H_2)$ with $v(H_1) \le n_1$, $v(H_2) \le n_2$ and for any $G_2 \subset G$ with $v(G_2) = v(H_2)$, there exists a subhypergraph $G_1 \subset G$ such that $G_1$ is an exact $(H_1, H_2)$-extension of $G_2$
and $(G_1, G_2)$ is $(K_1, K_2)$-maximal in $G$ for all $(K_1, K_2) \in \mathcal K^{\rho}$. 
\end{itemize}

Let $k \in \mathbb N$, $\rho > 0$, $1/\rho \in \left(s - 1 - \frac{1}{2^{k - s + 1}}, s - 1\right)$.
Set
$$
n_1(\rho) = \eta(\rho) + (k - s + 1) 2^{k - s + 1}, n_2(\rho) = \eta(\rho) + (k - s) 2^{k - s + 1}, n_3 = 2^{k - s + 1}, n_4 = 2.
$$

For a $s$-hypergraph $G$, denote by $d_G(x, \tilde G)$ the distance between a vertex $x \in V(G)$ and a subhypergraph $\tilde G \subset G$.
Denote by $d_G(\tilde G_1, \tilde G_2)$ the distance between $\tilde G_1, \tilde G_2 \subset G$.

\textbf{Proposition 3}.
\textit{
Let $G$ be $(n_1(\rho), n_2(\rho), n_3, n_4, \rho)$-sparse. Let $m \in \mathbb N$, $m \le 2^{k - s + 1}$. 
Then for any $\tilde G \subset G$ with $v(\tilde G) \le n_1(\rho)$ and any $x \in V(\tilde G)$, there exist a vertex $y \in V(G)$ and a path $P \subset G$ such that $d_G(y, \tilde G) = m$, $P$ connects $x$ and $y$, $e(P) = m$, $v(P) = e(P)(s - 1) + 1$.
Moreover, $G$ does not contain any path $P'$ of length at most $m$ such that $P'$ connects $y$ with a vertex from $\tilde G$ and $P' \neq P$.
} \\

\textbf{Proposition 4}.
\textit{
Let $G, H$ be $(n_1(\rho), n_2(\rho), n_3, n_4, \rho)$-sparse, where
$1/\rho \in \left(s - 1 - \frac{1}{2^{k - s + 1}}, s - 1\right) \setminus \mathcal Q_k$, 
$\mathcal Q_k = \{s - 1 - \frac{1}{2^{k - s + 1} + \frac{a}{b}} | a, b \in \mathbb N, a \le 2^{k - s + 1}\}$.
Then Duplicator has a winning strategy in $EHR(G, H, k)$.}

\subsection{Ehrenfeucht game}
\indent

Let us define the Ehrenfeucht game $\EHR(G_1,G_2;k)$ for hypergraphs $G_1, G_2$ and a number $k$. The game $\EHR(G_1,G_2;k)$ is a two-player game (the players are called Spoiler and Duplicator), which is played on a pair of nonempty disjoint $s$-hypergraphs $G_1$ and $G_2$ for a fixed number of rounds $k$. On the $i$-th round ($i=1,2,\ldots,k$) Spoiler selects an arbitrary vertex $x_i\in V(G_1)$ or $y_i\in V(G_2)$. After that Duplicator must select a vertex from the other graph. At the end of the $k$-th round there are $k$ chosen vertices in each of the hypergraphs: $x_1,x_2,\ldots,x_{k}\in V(G_1)$ and $y_1,y_2,\ldots,y_{k}\in V(G_2)$. Duplicator wins if and only if 
$$
G_1|_{\{x_1,x_2,\ldots,x_k\}} \cong G_2|_{\{y_1,y_2,\ldots,y_k\}}, 
$$
that is the induced subgraphs in the sets of chosen vertices are isomorphic.

The following theorem relates Ehrenfeucht game to the zero-one law. It is a corollary of Ehrenfeucht's Theorem (see~\cite{ehren}).

\textbf{Theorem 15}.
\textit{
The random $s$-hypergraph $G^{s}(n,p(n))$ obeys the zero-one $k$-law if and only if  
$$
 \lim\limits_{n,m\rightarrow\infty}\Pr[\text{Duplicator wins } \EHR(G_1,G_2;k)]=1,
$$
where $G_1\sim G^{s}(n,p(n))$, $G_2\sim G^{s}(m,p(m))$ are independently chosen and have disjoint vertex sets.
}

Thus, to prove that $G^{s}(n,p(n))$ obeys the zero-one law it is sufficient to describe an a.a.s. Duplicator's winning strategy in an Ehrenfeucht game with arbitrary but fixed number of rounds $k$.

\section{Proofs of theorems}
\indent

\textbf{Proof of theorem 6}.
Denote by $\mathcal B_{m, l}$ the set of all $l$-element subsets of the set $\{1, \ldots, m\}$.
Denote by $deg_G(x)$ the degree of a vertex $x \in V(G)$ in a hypergraph $G$.

By theorem 15 it is sufficient to describe Duplicator's a.a.s. winning strategy in $\EHR(G_1,G_2,k)$, where $G_1\sim G^{s}(n,n^{-\alpha})$, $G_2\sim G^{s}(m,m^{-\alpha})$ are independently chosen.

Let us consider several cases depending on the value of $\alpha$. In each case we find sets $\tilde \Omega_n \subset \Omega_n$ such that $\Pr[G^s(n, n^{-\alpha}) \in \tilde \Omega_n] \to 1$, $n \to \infty$,
and prove that for any $n_1, n_2 \in \mathbb N$ and hypergraphs $G \in \tilde \Omega_{n_1}$, $H \in \tilde \Omega_{n_2}$ Duplicator has a winning strategy in the game $EHR(G, H, k)$. Denote the vertices chosen in $G$ and $H$ in the $i$-th round by $x_i$ and $y_i$ respectively.

1) Let $1/\alpha > {{k - 1} \choose {s - 1}}$.
Let $\tilde \Omega_n$ be the set of all hypergraphs $G \in \Omega_n$ satisfying the following property: for any $r \in \{s - 1, \ldots, k - 1\}$, $\mathcal A \subset \mathcal B_{r, s - 1}$  and distinct vertices $z_1, \ldots, z_r \in V_n$, there exists a vertex $z \in V_n \setminus \{z_1, \ldots, z_{r}\}$ such that $\{z_{i_1}, \ldots, z_{i_{s - 1}}, z\} \in E(G)$ for all $\{i_1, \ldots, i_{s - 1}\} \in \mathcal A$ and $\{z_{i_1}, \ldots, z_{i_{s - 1}}, z\} \notin E(G)$ for all $\{i_1, \ldots, i_{s - 1}\} \in \mathcal B_{r, s - 1} \setminus  \mathcal A$. This property is called the \textit{full level $(k - 1)$ extension property}. By theorem 14 
$\lim_{n \to \infty} \Pr[G^s(n, n^{-\alpha}) \in \tilde \Omega_n] = 1$. 

Let us describe the winning strategy of Duplicator in $EHR(G, H, k)$, where $G \in \tilde \Omega_{n_1}$, $H \in \tilde \Omega_{n_2}$.
Let Spoiler choose a vertex $x_i \in V(G)$ in the $i$-th round, $1 \le i \le k$.
Then by full level $(i - 1)$ extension property Duplicator can choose a vertex $y_i$ such that for any distinct $j_1, \ldots, j_{s  - 1} \in \{1, \ldots, i - 1\}$ we have
$\{x_{j_1}, \ldots, x_{j_{s - 1}}, x_i\} \in E(G)$ if and only if $\{y_{j_1}, \ldots, y_{j_{s - 1}}, y_i\} \in E(H)$. Therefore, Duplicator wins after the $k$-th round.

2) Let ${{k - 1} \choose {s - 1}} - 1 < 1/\alpha \le {{k - 1} \choose {s - 1}}$.
Let $\tilde \Omega_n$ be the set of all hypergraphs $G \in \Omega_n$ satisfying the following properties:
\begin{itemize}
\item[1)] $G$ satisfies full level $(k - 2)$ extension property,
\item[2)] for any distinct vertices $z_1, \ldots, z_{k - 1} \in V_n$ and $\mathcal A \subsetneq \mathcal B_{k - 1, s - 1}$,  there exists a vertex $z \in V_n \setminus \{z_1, \ldots, z_{k - 1}\}$ such that $\{z_{i_1}, \ldots, z_{i_{s - 1}}, z\} \in E(G)$ for all $\{i_1, \ldots, i_{s - 1}\} \in \mathcal A$ and $\{z_{i_1}, \ldots, z_{i_{s - 1}},$ $z\} \notin E(G)$ for all $\{i_1, \ldots, i_{s - 1}\} \in \mathcal B_{k - 1, s - 1} \setminus \mathcal A$,
\item[3)] for any distinct vertices $z_1, \ldots, z_{k - 2}$ and $\mathcal A \subset \mathcal B_{k - 2, s - 1}$ there exist distinct vertices $z_{k - 1}, z_k \in V_n \setminus \{z_1, \ldots, z_{k - 2}\}$ such that 
$\{z_{i_1}, \ldots, z_{i_{s - 1}}, z_{k - 1}\} \in E(G)$ for all $\{i_1, \ldots, i_{s - 1}\} \in \mathcal A$, $\{z_{i_1}, \ldots,$ $z_{i_{s - 1}}, z_{k - 1}\} \notin E(G)$ for all $\{i_1, \ldots, i_{s - 1}\} \in \mathcal B_{k - 2, s - 1} \setminus \mathcal A$ and
$\{z_{i_1}, \ldots, z_{i_{s - 1}}, z_k\} \in E(G)$ for all $\{i_1, \ldots, i_{s - 1}\} \in \mathcal B_{k - 1, s - 1}$.
\item[4)] for any distinct vertices $z_1, \ldots, z_{k - 2}$ and $\mathcal A \subset \mathcal B_{k - 2, s - 1}$ there exists a vertex $z_{k - 1} \in V_n \setminus \{z_1, \ldots, z_{k - 2}\}$ such that 
$\{z_{i_1}, \ldots, z_{i_{s - 1}}, z_{k - 1}\} \in E(G)$ for all $\{i_1, \ldots, i_{s - 1}\} \in \mathcal A$, $\{z_{i_1}, \ldots,$ $z_{i_{s - 1}}, z_{k - 1}\} \notin E(G)$ for all $\{i_1, \ldots, i_{s - 1}\} \in \mathcal B_{k - 2, s - 1} \setminus \mathcal A$ and
there does not exist a vertex $z_k \in V_n \setminus \{z_1, \ldots, z_{k - 1}\}$ such that $\{z_{i_1}, \ldots, z_{i_{s - 1}}, z_k\} \in E(G)$ for all $\{i_1, \ldots, i_{s - 1}\} \in \mathcal B_{k - 1, s - 1}$.
\end{itemize}

Note that properties 1)-3) mean that the hypergraph $G$ satisfies $L^0_{(\tilde G, \tilde H)}$ for some pairs $(\tilde G, \tilde H)$, where $\rho^{\max}(\tilde G, \tilde H)$ equals ${{k - 2} \choose {s - 1}}$,
${{k - 1} \choose {s - 1}} - 1$ and $\frac{|\mathcal A| + {{k - 1} \choose {s - 1}}}{2}$, where $|\mathcal A| \le {{k - 2} \choose {s - 1}}$, respectively. Since all these values don't exceed ${{k - 1} \choose {s - 1}} - 1 < 1/\alpha$, theorem 14 implies that the random hypergraph
$G^s(n, n^{-\alpha})$ satisfies properties 1)-3) a.a.s. Since $1/\alpha \le {{k - 1} \choose {s - 1}}$, the property 4) is satisfied if the hypergraph $G$ satisfies $L^1_{(\tilde G, \tilde H)}$ for some pair
$(\tilde G, \tilde H)$, where $\rho^{\max}(\tilde G, \tilde H) = |\mathcal A| \le {{k - 2} \choose {s - 1}} < 1/\alpha$.
Therefore, $\lim_{n \to \infty} \Pr[G^s(n, n^{-\alpha}) \in \tilde \Omega_n] = 1$.

Let us describe the winning strategy of Duplicator in $EHR(G, H, k)$, where $G \in \tilde \Omega_{n_1}$, $H \in \tilde \Omega_{n_2}$.
Let Spoiler choose a vertex $x_i \in V(G)$ in the $i$-th round, $1 \le i \le k - 2$. 
Then by full level $(i - 1)$ extension property Duplicator chooses a vertex $y_i \in V(H)$ such that for any distinct $j_1, \ldots, j_{s  - 1} \in \{1, \ldots, i - 1\}$ we have
$\{x_{j_1}, \ldots, x_{j_{s - 1}}, x_i\} \in E(G)$ if and only if $\{y_{j_1}, \ldots, y_{j_{s - 1}}, y_i\} \in E(H)$.

Let Spoiler choose a vertex $x_{k - 1} \in V(G)$ in the $(k - 1)$-th round. Since $H$ satisfies 3) and 4) there exists a vertex $y_{k - 1} \in V(H)$ such that $H|_{\{y_1, \ldots, y_{k - 1}\}}$ and $G|_{\{x_1, \ldots, x_{k - 1}\}}$ are isomorphic
and there exists a vertex $x_k \in V(G)$ such that $deg_{G|_{\{x_1, \ldots, x_k\}}}(x_k) = {{k - 1} \choose {s - 1}}$
if and only if there exists a vertex $y_k \in V(H)$ such that $deg_{H|_{\{y_1, \ldots, y_k\}}}(y_k) = {{k - 1} \choose {s - 1}}$.

Let Spoiler choose a vertex $x_{k} \in V(G)$ in the $k$-th round. If $deg_{G|_{\{x_1, \ldots, x_k\}}}(x_k)  = {{k - 1} \choose {s - 1}}$, then by the definition
of the vertex chosen by Duplicator in the $(k - 1)$-th round there exists a vertex $y_k \in V(H)$ with $deg_{H|_{\{y_1, \ldots, y_k\}}}(y_k) = {{k - 1} \choose {s - 1}}$. If $deg_{G|_{\{x_1, \ldots, x_k\}}}(x_k)  < {{k - 1} \choose {s - 1}}$, then by property 2) there exists a vertex $y_k \in V(H)$
such that for any distinct $j_1, \ldots, j_{s  - 1} \in \{1, \ldots, i - 1\}$ we have
$\{x_{j_1}, \ldots, x_{j_{s - 1}}, x_i\} \in E(G)$ if and only if $\{y_{j_1}, \ldots, y_{j_{s - 1}}, y_i\} \in E(H)$. So in both cases Duplicator wins.

3) Let  ${{k - 1} \choose {s - 1}} - 1 - \frac{s - 1}{k - 1} + \frac{2\left(1 + \frac{s - 1}{k - 1}\right)}{{{k - 1} \choose {s - 1}} + 2} < 1/\alpha \le {{k - 1} \choose {s - 1}} - 1$.

Let $z_1, \ldots, z_k \in V_n$.
Let $L_C(z_1, \ldots, z_k)$, where $C \in \mathcal B_{k - 1, s - 1}$, be the following property of a hypergraph $G \in \Omega_n$: 
$\{z_{i_1}, \ldots, z_{i_{s - 1}}, z_k\} \in E(G)$ for all $\{i_1, \ldots, i_{s - 1}\} \in \mathcal B_{k - 1, s - 1} \setminus \{C\}$
and $\{z_{i_1}, \ldots, z_{i_{s - 1}}, z_k\} \notin E(G)$ for $\{i_1, \ldots, i_{s - 1}\} = C$.
Let $L_{\varnothing}(z_1, \ldots, z_k)$ be the following property of a hypergraph $G \in \Omega_n$: 
$\{z_{i_1}, \ldots, z_{i_{s - 1}}, z_k\} \in E(G)$ for all $\{i_1, \ldots, i_{s - 1}\} \in \mathcal B_{k - 1, s - 1}$.

Let $\tilde \Omega_n$ be the set of all hypergraphs $G \in \Omega_n$ satisfying the following properties:
\begin{itemize}
\item[1)] $G$ satisfies full level $(k - 2)$ extension property,
\item[2)] for any distinct vertices $z_1, \ldots, z_{k - 1} \in V_n$ and $\mathcal A \subset \mathcal B_{k - 1, s - 1}$, where $|\mathcal A| \le {{k - 1} \choose {s - 1}} - 2$, there exists a vertex $z \in V_n \setminus \{z_1, \ldots, z_{k - 1}\}$ such that $\{z_{i_1}, \ldots, z_{i_{s - 1}}, z\} \in E(G)$ for all $\{i_1, \ldots, i_{s - 1}\} \in \mathcal A$ and $\{z_{i_1}, \ldots, z_{i_{s - 1}}, z\} \notin E(G)$ for all $\{i_1, \ldots, i_{s - 1}\} \in \mathcal B_{k - 1, s - 1} \setminus \mathcal A$,
\item[3)] for any distinct vertices $z_1, \ldots, z_{k - 2}$, $\mathcal A \subset \mathcal B_{k - 2, s - 1}$ and $\mathcal C \subset \mathcal B_{k - 1, s - 1} \cup \{\varnothing\}$ there exists a vertex $z_{k - 1} \in V_n \setminus \{z_1, \ldots, z_{k - 2}\}$ such that  
$\{z_{i_1}, \ldots, z_{i_{s - 1}}, z_{k - 1}\} \in E(G)$ for all $\{i_1, \ldots, i_{s - 1}\} \in \mathcal A$, $\{z_{i_1}, \ldots, z_{i_{s - 1}}, z_{k - 1}\} \notin E(G)$ for all $\{i_1, \ldots, i_{s - 1}\} \in \mathcal B_{k - 2, s - 1} \setminus \mathcal A$ and for any $C \in \mathcal C$ there exists a vertex $z_k^C \in V_n \setminus \{z_1, \ldots, z_{k - 1}\}$ such that the property 
$L_C(z_1, \ldots, z_{k - 1}, z_k^C)$ is satisfied and for any $C \in (\mathcal B_{k - 1, s - 1} \cup \{\varnothing\}) \setminus \mathcal C$ there does not exist
a vertex $z_k^C \in V_n \setminus \{z_1, \ldots, z_{k - 1}\}$ such that $L_C(z_1, \ldots, z_{k - 1}, z_k^C)$ is satisfied.
\end{itemize}

Note that properties 1)-2) mean that the hypergraph $G$ satisfies $L^0_{(\tilde G, \tilde H)}$ for some pairs $(\tilde G, \tilde H)$, where $\rho^{\max}(\tilde G, \tilde H)$ equals ${{k - 2} \choose {s - 1}}$ and ${{k - 1} \choose {s - 1}} - 2$ respectively.  Since $1/\alpha \le {{k - 1} \choose {s - 1}} - 1$, the property 4) is satisfied if the hypergraph $G$ satisfies $L^{1}_{(\tilde G, \tilde H)}$ for some pair
$(\tilde G, \tilde H)$, where 
\begin{multline*}
\rho^{\max}(\tilde G, \tilde H) \le \frac{{{k - 2} \choose {s - 1}} + {{k - 1} \choose {s - 1}} + {{k - 1} \choose {s - 1}} \left({{k - 1} \choose {s - 1}} - 1\right)}
{{{k - 1} \choose {s - 1}} + 2} = {{k - 1} \choose {s - 1}} - 1 - \frac{{{k - 2} \choose {s - 2}} - 2}{{{k - 1} \choose {s - 1}} + 2} = \\
= {{k - 1} \choose {s - 1}} - 1 - \frac{s - 1}{k - 1} + \frac{2\left(1 + \frac{s - 1}{k - 1}\right)}{{{k - 1} \choose {s - 1}} + 2}
\end{multline*}
Therefore, $\lim_{n \to \infty} \Pr[G^s(n, n^{-\alpha}) \in \tilde \Omega_n] = 1$. 

Let us describe the winning strategy of Duplicator in $EHR(G, H, k)$, where $G \in \tilde \Omega_{n_1}$, $H \in \tilde \Omega_{n_2}$.
Let Spoiler choose a vertex $x_i \in V(G)$ in the $i$-th round, $1 \le i \le k - 2$. 
Then by full level $(i - 1)$ extension property Duplicator chooses a vertex $y_i \in V(H)$ such that for any distinct $j_1, \ldots, j_{s  - 1} \in \{1, \ldots, i - 1\}$ we have
$\{x_{j_1}, \ldots, x_{j_{s - 1}}, x_i\} \in E(G)$ if and only if $\{y_{j_1}, \ldots, y_{j_{s - 1}}, y_i\} \in E(H)$.

Let Spoiler choose a vertex $x_{k - 1} \in V(G)$ in the $(k - 1)$-th round. Then by property 3) there exists a vertex $y_{k - 1} \in V(H)$ such that
$H|_{\{y_1, \ldots, y_{k - 1}\}}$ and $G|_{\{x_1, \ldots, x_{k - 1}\}}$ are isomorphic and for each $C \in \mathcal B_{k - 1, s - 1} \cup \{\varnothing\}$
there exists a vertex $x_k \in V(G)$ such that $G$ satisfies $L_C(x_1, \ldots, x_k)$ if and only if there exists a vertex 
$y_k \in V(G)$ such that $H$ satisfies $L_C(y_1, \ldots, y_k)$.

Let Spoiler choose a vertex $x_{k} \in V(G)$ in the $k$-th round. If $deg_{G|_{\{x_1, \ldots, x_k\}}}(x_k)  \ge {{k - 1} \choose {s - 1}} - 1$, then there exists exactly one $C \in \mathcal B_{k - 1, s - 1} \cup \{\varnothing\}$ such that $G$ satisfies $L_C(x_1, \ldots, x_k)$. By the definition
of the vertex chosen by Duplicator in the $(k - 1)$-th round there exists a vertex $y_k \in V(H)$ such that $H$ satisfies $L_C(y_1, \ldots, y_k)$. 
If $deg_{G|_{\{x_1, \ldots, x_k\}}}(x_k)  < {{k - 1} \choose {s - 1}} - 1$, then by property 2) there exists a vertex $y_k \in V(H)$
such that for any distinct $j_1, \ldots, j_{s  - 1} \in \{1, \ldots, i - 1\}$ we have
$\{x_{j_1}, \ldots, x_{j_{s - 1}}, x_i\} \in E(G)$ if and only if $\{y_{j_1}, \ldots, y_{j_{s - 1}}, y_i\} \in E(H)$. So in both cases Duplicator wins. \\

\textbf{Proof of theorem 7}.
Denote by $\mathcal B_{m, l}$ the set of all $l$-element subsets of the set $\{1, \ldots, m\}$.
Let $\mathcal A = \{\{A, B\}: A, B \in \mathcal B_{k - 2, s - 1}, A \neq B\}$.

Let $L$ be a property expressed by the formula 
$$
\exists x_1 \ldots \exists x_{k - 2} \bigwedge_{\{A, B\} \in \mathcal A} \left(\bigwedge_{C \in \mathcal B_{k - 1, s - 1}} \exists x_{k - 1} \mbox{ } Q_{A, B, C}(x_1, \ldots, x_{k - 1})\right),
$$
where
$$Q_{A, B, C}(x_1, \ldots, x_{k - 1}) = R^2_{A, B}(x_1, \ldots, x_{k - 1}) \land W_C(x_1, \ldots, x_{k - 1}),$$
\begin{multline*}
R^2_{A, B}(x_1, \ldots, x_{k - 1}) = \left(\bigwedge_{\{i_1, \ldots, i_{s - 1}\} \in \mathcal B_{k - 2, s - 1} \setminus \{A, B\}} N(x_{i_1}, \ldots, x_{i_{s - 1}}, x_{k - 1}) \right) \bigwedge \\
\bigwedge  \left(\bigwedge_{\{i_1, \ldots, i_{s - 1}\} \in \{A, B\}} (\neg N(x_{i_1}, \ldots, x_{i_{s - 1}}, x_{k - 1}))\right) \bigwedge
((x_1 \neq x_2) \land \ldots \land (x_{k - 2} \neq x_{k - 1})),
\end{multline*}
$$
W_C(x_1, \ldots, x_{k - 1}) = \left(\bigwedge_{C' \in \mathcal B_{k - 1, s - 1} \setminus \{C\}} \exists x_k \mbox{ } R^1_{C'}(x_1, \ldots, x_k)\right) \land 
(\neg (\exists x_k \mbox{ } R^1_{C}(x_1, \ldots, x_k)),
$$
\begin{multline*}
R^1_{C}(x_1, \ldots, x_k) = \left(\bigwedge_{\{i_1, \ldots, i_{s - 1}\} \in \mathcal B_{k - 1, s - 1} \setminus \{C\}} N(x_{i_1}, \ldots, x_{i_{s - 1}}, x_{k})\right) \bigwedge \\ 
\bigwedge (\neg N(x_{j_1}, \ldots, x_{j_{s - 1}}, x_k)) \bigwedge ((x_ k \neq x_1) \land \ldots \land (x_k \neq x_{k - 1})),
\end{multline*}
where $\{j_1, \ldots, j_{s - 1}\} = C$.

Let $K$ be a hypergraph with
\begin{multline*}
V(K) = \{x_1, \ldots, x_{k - 2}\} \cup \left\{x_{k - 1}^{A, B, C}: \{A, B\} \in \mathcal A, C \in \mathcal B_{k - 1, s - 1} \right\} \cup \\
\cup \left\{x_k^{A, B, C, C'}: \{A, B\} \in \mathcal A, C, C' \in \mathcal B_{k - 1, s - 1}, C \neq C' \right\}
\end{multline*}
and 
\begin{multline*}
E(K) = \left\{\left\{x_{i_1}, \ldots, x_{i_{s - 1}}, x_{k - 1}^{A, B, C}\right\}: \{A, B\} \in \mathcal A, C \in \mathcal B_{k - 1, s - 1}, \{i_1, \ldots, i_{s - 1}\} \in \mathcal B_{k - 2, s - 1} \setminus \{A, B\}\right\} \cup \\
\cup \Bigl\{\left\{x_{i_1}, \ldots, x_{i_{s - 1}}, x_k^{A, B, C, C'}\right\}: \{A, B\} \in \mathcal A, C, C' \in \mathcal B_{k - 1, s - 1}, C \neq C', 
\\
\{i_1, \ldots, i_{s - 1}\} \in \mathcal B_{k - 2, s - 1} \setminus \{C'\} \Bigr\}
\cup \Bigl\{\left\{x_{i_1}, \ldots, x_{i_{s - 2}}, x^{A, B, C}_{k - 1}, x_k^{A, B, C, C'}\right\}: \\
\{A, B\} \in \mathcal A, C, C' \in \mathcal B_{k - 1, s - 1}, C \neq C', \{i_1, \ldots, i_{s - 2}\} \in \mathcal B_{k - 2, s - 2},
\{i_1, \ldots, i_{s - 2}, k - 1\} \neq C'\Bigr\}.
\end{multline*}

The hypergraph $K$ is strictly balanced and 
$$|V(K)| = k - 2 + |\mathcal A| {{k - 1} \choose {s - 1}} + |\mathcal A| {{k - 1} \choose {s - 1}} \left({{k - 1} \choose {s - 1}} - 1\right),$$
$$|E(K)| = |\mathcal A| {{k - 1} \choose {s - 1}} \left({{k - 2} \choose {s - 1}} - 2\right) + |\mathcal A| {{k - 1} \choose {s - 1}} \left({{k - 1} \choose {s - 1}} - 1\right)^2.$$

Let us prove that if a hypergraph $G$ posseses property $L$, then it contains either a subhypergraph isomorphic to $K$ or a
subhypergraph $\tilde K$ with $v(\tilde K) \le v(K)$ and $\rho(\tilde K) > \rho(K)$. 
Suppose that a hypergraph $G$ posseses property $L$. Then there exist
distinct vertices $x_1, \ldots, x_{k - 2}$, $x_{k - 1}^{A, B, C} \in V(G)$, where 
$\{A, B\} \in \mathcal A$, $C \in \mathcal B_{k - 1, s - 1}$, such that the predicates $Q_{A, B, C}(x_1, \ldots, x_{k - 2}, x_{k - 1}^{A, B, C})$ are true for all
$\{A, B\} \in \mathcal A$, $C \in \mathcal B_{k - 1, s - 1}$. Let $\tilde V =  \{x_1, \ldots, x_{k - 2}\} \cup \left\{x_{k - 1}^{A, B, C}: \{A, B\} \in \mathcal A, C \in \mathcal B_{k - 1, s - 1} \right\}$.
The truth of the predicate $R^2_{A, B}(x_1, \ldots,$ $x_{k - 2}, x_{k - 1}^{A, B, C})$ implies
that $deg_{G|_{\left\{x_1, \ldots, x_{k - 2}, x_{k - 1}^{A, B, C}\right\}}}(x_{k - 1}^{A, B, C}) = {{k - 2} \choose {s - 1}} - 2$.
The truth of the predicates $W_C(x_1, \ldots, x_{k - 1}^{A, B, C})$ implies that there exist distinct vertices 
$x_k^1, \ldots, x_k^{N_1}, \tilde x_k^1, \ldots, \tilde x_k^{N_2} \in V(G) \setminus \tilde V$ and numbers 
$v_1, \ldots, v_{N_1}$, $w_1, \ldots, w_{N_2} \in \mathbb N$ such that 
\newline
$deg_{G|_{\tilde V \cup \{x_k^i\}}}(x_k^i) = {{k - 2} \choose {s - 1}} - 1 + v_i {{k - 2}\choose {s - 2}}$
and
$deg_{G|_{\tilde V \cup \{\tilde x_k^i\}}}(\tilde x_k^i) = {{k - 2} \choose {s - 1}} + w_i \left({{k - 2}\choose {s - 2}} - 1 \right)$,
where 
$$\sum_{i = 1}^{N_1} v_i = |\mathcal A| {{k - 2}\choose {s - 2}} {{k - 2}\choose {s - 1}} + |\mathcal A| {{k - 2}\choose {s - 1}} \left({{k - 2}\choose {s - 1}} - 1\right) = |\mathcal A| {{k - 2}\choose {s - 1}} \left({{k - 1}\choose {s - 1}} - 1\right),$$ 
$$\sum_{i = 1}^{N_1} w_i  =|\mathcal A| {{k - 2}\choose {s - 2}} \left({{k - 2}\choose {s - 2}} - 1\right) + |\mathcal A| {{k - 2}\choose {s - 1}} {{k - 2}\choose {s - 2}} = |\mathcal A| {{k - 2}\choose {s - 2}} \left({{k - 1}\choose {s - 1}} - 1\right).$$

Let $\tilde K = G|_{\tilde V \cup \{x_k^1, \ldots, x_k^{N_1}, \tilde x_k^1, \ldots, \tilde x_k^{N_2}\}}$. Then
$$V(\tilde K) = k - 2 + |\mathcal A| {{k - 1}\choose {s - 1}} + N_1 + N_2,$$
\begin{multline*}
|E(\tilde K)| \ge |\mathcal A| {{k - 1} \choose {s - 1}} \left({{k - 2} \choose {s - 1}} - 2\right) + N_1 \left({{k - 2} \choose {s - 1}} - 1\right) + N_2 {{k - 2} \choose {s - 1}} +
\sum_{i = 1}^{N_1} v_i {{k - 2} \choose {s - 2}} + \\
+ \sum_{i = 1}^{N_2} w_i \left({{k - 2} \choose {s - 2}} - 1\right) = 
 |\mathcal A| {{k - 2}\choose {s - 2}} {{k - 1}\choose {s - 1}} \left({{k - 1}\choose {s - 1}} - 1 \right) + \\ +
|\mathcal A| {{k - 1}\choose {s - 1}} \left({{k - 2} \choose {s - 1}} - 2\right) - |\mathcal A| {{k - 2}\choose {s - 2}} \left({{k - 1}\choose {s - 1}} - 1\right) +
N_1 \left({{k - 2} \choose {s - 1}} - 1\right) + N_2 {{k - 2} \choose {s - 1}}.
\end{multline*}
Hence
$$
\rho(\tilde K) \ge \frac{e_0 + N_1 \left({{k - 2} \choose {s - 1}} - 1\right) + N_2 {{k - 2} \choose {s - 1}}}
{k - 2 + |\mathcal A| {{k - 1}\choose {s - 1}} + N_1 + N_2},
$$
where 
$$e_0 = |\mathcal A| {{k - 2}\choose {s - 2}} {{k - 1}\choose {s - 1}} \left({{k - 1}\choose {s - 1}} - 1 \right) + 
|\mathcal A| {{k - 1}\choose {s - 1}} \left({{k - 2} \choose {s - 1}} - 2\right) - |\mathcal A| {{k - 2}\choose {s - 2}} \left({{k - 1}\choose {s - 1}} - 1\right).$$
Note that $1 \le N_1 \le N_1^{\max}$, $1 \le N_2 \le N_2^{\max}$, where 
$N_1^{\max} = |\mathcal A| {{k - 2}\choose {s - 1}} \left({{k - 1}\choose {s - 1}} - 1\right)$, $N_2^{\max} = |\mathcal A| {{k - 2}\choose {s - 2}} \left({{k - 1}\choose {s - 1}} - 1\right)$.
Since
$$
 \frac{e_0 + N_1 \left({{k - 2} \choose {s - 1}} - 1\right)}
{k - 2 + |\mathcal A| {{k - 1}\choose {s - 1}} + N_1} > {{k - 2} \choose {s - 1}},
$$
the density $\rho(\tilde K)$ decreases as $N_2$ grows. Similarly, $\rho(\tilde K)$ decreases as $N_1$ grows.
Therefore,
\begin{multline*}
\rho(\tilde K) \ge 
\frac{e_0 + (N_1^{\max} + N_2^{\max}){{k - 2} \choose {s - 1}} - N_1^{\max}}{k - 2 + |\mathcal A| {{k - 1} \choose {s - 1}} + N_1^{\max} + N_2^{\max}} = \\
=
\frac{ |\mathcal A| {{k - 1}\choose {s - 1}}^2 \left({{k - 1}\choose {s - 1}} - 1 \right) +
|\mathcal A| {{k - 1}\choose {s - 1}} \left({{k - 2} \choose {s - 1}} - 2\right) - |\mathcal A| {{k - 1}\choose {s - 1}} \left({{k - 1}\choose {s - 1}} - 1\right)}{ k - 2 + |\mathcal A| {{k - 1} \choose {s - 1}}^2} 
= \rho(K).
\end{multline*}
The density $\rho(\tilde K)$ equals $\rho(K)$ if and only if $N_1 = N_1^{\max}$, $N_2 = N_2^{\max}$ and $\tilde K$ is isomorphic to $K$.

We have
\begin{multline*}
\rho(K) = {{k -1} \choose {s - 1}} - 1 - \frac{|\mathcal A| {{k - 1}\choose {s - 1}} \left({{k - 2}\choose {s - 2}} + 1\right) + (k - 2) \left({{ k -1} \choose {s - 1}} - 1\right)}{k - 2 + |\mathcal A| {{k - 1} \choose {s - 1}}^2} \ge \\ \ge
 {{ k -1} \choose {s - 1}} - 1 - \frac{{{k - 2}\choose {s - 2}} + 1}{{{k - 1} \choose {s - 1}}} - \frac{(k - 2)}{|\mathcal A| {{k - 1} \choose {s - 1}}} \ge
 {{ k -1} \choose {s - 1}} - 1 - \frac{s - 1}{k - 1} - \frac{2}{{{k - 1} \choose {s - 1}}}.
\end{multline*}

Let $\alpha = 1/\rho(K)$. Then by theorem 12 $\Pr[G^s(n, n^{-\alpha}) \models L] = \Pr[G^s(n, n^{-\alpha}) \models L_K] + o(1)$.
By theorem 13 there exists $\lim_{n \to \infty} \Pr[G^s(n, n^{-\alpha}) \models L_K] \in (0, 1)$.
Hence there exists 
\newline
$\lim_{n \to \infty} \Pr[G^s(n, n^{-\alpha}) \models L] \in (0, 1)$.
So $G^s(n, n^{-\alpha})$ does not follow the zero-one $k$-law and the theorem is proved. \\

\textbf{Proof of Theorem 8}.
Theorems 12 and 14 imply that $G^s(n, n^{-\alpha})$  is $(n_1(\rho), n_2(\rho), n_3, n_4)$-sparse a.a.s. 
Then by Proposition 4 Duplicator has a winning strategy in the game $EHR(G^s(n_1, n_1^{-\alpha}),$ $G^s(n_2, n_2^{-\alpha}), k)$ a.a.s.
Therefore, by theorem 15 $G^s(n, n^{-\alpha})$ follows the zero-one $k$-law. \\

\textbf{Proof of Theorem 9}.
Let 
$$D_1(x_1, x_2) = (x_1 = x_2) \lor (\exists x_3 \ldots \exists x_s \quad N(x_1, \ldots, x_s)),$$ 
$$D_{i}(x_1, x_2) = \exists x_3 \mbox{ } (D_{\lfloor i/2 \rfloor}(x_1, x_3) \land D_{\lceil i/2 \rceil}(x_3, x_2)), \quad i > 1.$$ 
Let $$D^{=}_1(x_1, x_2) = D_1(x_1, x_2) \land (x_1 \neq x_2), \quad D^{=}_i(x_1, x_2) = D_i(x_1, x_2) \land (\neg (D_{i - 1}(x_1, x_2))), \quad i > 1.$$ 
Then quantifier depths of formulae $D_i$ and $D^{=}_i$ equal $\lceil \log_2 i \rceil + s - 2$.
$D_i(x_1, x_2)$ expresses the property that the the distance between vertices $x_1$ and $x_2$ is at most $i$. $D^{=}_i(x_1, x_2)$ expresses the property that the distance between $x_1$ and $x_2$ equals $i$.
Let
\begin{multline*}
\tilde D_1(x, x_1, x_2) = (x_1 \neq x) \land (x_2 \neq x) \land ((x_1 = x_2) \lor (\exists x_3 \ldots \exists x_s \quad (N(x_1, \ldots, x_s) \land \\ \land
(x_3 \neq x) \land \ldots \land (x_s \neq x)))),
\end{multline*}
$$
\tilde D_i(x, x_1, x_2) =  (x_1 \neq x) \land (x_2 \neq x) \land (\exists x_3 \mbox{ } ((x_3 \neq x) \land \tilde D_{\lfloor i/2 \rfloor}(x, x_1, x_3) \land \tilde D_{\lceil i/2 \rceil}(x, x_3, x_2))), \quad i > 1.
$$
Then $\tilde D_i(x, x_1, x_2)$ expresses the property that there exists a path of length at most $i$ which connects $x_1$ and $x_2$ and does not contain $x$. 
The quantifier depth of $\tilde D_i$  equals $\lceil \log_2 i \rceil + s - 2$.

For $i \in \mathbb N$, let
$$
C_i(x_1) = \exists x_2 \mbox{ } (D^{=}_{i}(x_1, x_2) \land 
(\exists x_3 \mbox{ } (B_{i + 1}(x_1, x_2, x_3) \land \tilde D_i(x_3, x_1, x_2)))),
$$
where 
$$
B_i(x_1, x_2, x_3) =  D^{=}_{\lfloor i/2 \rfloor}(x_1, x_3) \land D^{=}_{\lceil i/2 \rceil}(x_3, x_2), \quad i \ge 2.
$$

Then the quantifier depth of $C_i$ equals $\lceil \log_2 i \rceil + s$.

Let $a \in \mathbb N$, $a \le 2^{k - s - 2} + 2^{k - s - 3} + 1$. 
Then there exist $a_1, a_2, a_3 \in \mathbb N$ such that $2 a_1 + 2a_2 + 1 + a_3 = 2^{k - s + 1} + a$,
$2 \le a_1 \le 2^{k - s}$, $a_2 \le 2^{k - s - 4}$, $a_3 \le 2^{k - s - 2}$ and $a_2 < a_1$.

Let $L$ be a property expressed by the formula 
$$
\exists x_1 \exists x_2 \mbox{ } (D^{=}_{a_1}(x_1, x_2) \land (\exists x_3 \mbox{ } (B_{a_1}(x_1, x_2, x_3) \land Q(x_3))) \land
 (\exists x_3 \mbox{ } (B_{a_1}(x_1, x_2, x_3) \land (\neg Q(x_3)))),
$$
where
$$
Q(x_3) = \exists x_4 \mbox{ } (D^{=}_{a_3}(x_3, x_4) \land C_{a_2}(x_4)).
$$
Then the quantifier depth of the formula expressing the property $L$ is at most $k$.

Let us prove that $\Pr[G^s(n, n^{-\alpha}) \models L]$ tends neither to $0$ nor to $1$.

Let $K$ be a hypergraph with
$$V(K) = \{x_1, \ldots, x_{2a_1 (s - 1)}, y_1, \ldots, y_{(2a_2 + 1)(s - 1)}, z_{1, \ldots, z_{a_3 (s - 1)} - 1}\},$$
\begin{multline*}
E(K) = \{x_{(i - 1)(s - 1) + 1}, \ldots, x_{i(s - 1) + 1}: 1 \le i < 2a_1\} \cup \{\{x_{(2a_1 - 1)(s - 1) + 1}, \ldots, x_{2a_1 (s - 1)}, x_1\}\} \cup \\
\cup \{y_{(i - 1)(s - 1) + 1}, \ldots, y_{i(s - 1) + 1}: 1 \le i < 2a_2 + 1\} \cup \{\{y_{2a_2 (s - 1) + 1}, \ldots, y_{(2a_2 + 1) (s - 1)}, y_1\}\} \cup \\
\cup \{\{x_1, z_1, \ldots, z_{s - 1}\}, \{z_{(a_3 - 1)(s - 1)}, \ldots, z_{a_3(s - 1) - 1}, y_1\}\} \cup \{z_{i(s - 1)}, \ldots, z_{(i + 1)(s - 1)}: 1 \le i < a_3 - 1\}
\end{multline*}

The hypergraph $K$ is strictly balanced and $1/\rho(K) = s - 1 - \frac{1}{2^{k - s + 1} + a} = \alpha$.
Theorem 12 implies that $\Pr[G^s(n, n^{-\alpha}) \models L] \ge  \Pr[G^s(n, n^{-\alpha}) \models L_K] + o(1)$.
Since by theorem 13 there exists $\lim_{n \to \infty} \Pr[G^s(n, n^{-\alpha}) \models L_K] \in (0, 1)$, we obtain that $\lim\inf_{n \to \infty} \Pr[G^s(n, n^{-\alpha}) \models L] > 0$.

Let us show that if a hypergraph posseses property $L$, then it contains a subhypergraph $\tilde K$ with $e(\tilde K) \le 2^{k - s + 1} + a$ and $v(\tilde K) \le e(\tilde K)(s - 1) - 1$.
Let $G$ be a hypergraph which posseses property $L$. Then there exist distinct vertices $x_1, x_2, x_3^1, x_3^2 \in V(G)$ such that 
$$d_G(x_1, x_2) = a_1, \quad d_G(x_1, x_3^i) = \lfloor a_1/2 \rfloor, d_G(x_2, x_3^i) = \lceil a_1 / 2 \rceil, \quad i \in \{1, 2\},$$
the predicate $Q(x_3^1)$ is true and the predicate $Q(x_3^2)$ is false.
Let $W_1^i$ be a path of length $\lfloor a_1/2 \rfloor$ connecting $x_1$ and $x_3^i$, $W_2^i$ be a path of length $\lceil a_1 / 2 \rceil
$ connecting $x_2$ and $x_3^i$. Then $W_1^1 \cup W_2^1 \cup W_1^2 \cup W_2^2$ contains either a subhypergraph $\tilde K$ with
$e(\tilde K) \le 2 a_1$ and $v(\tilde K) \le e(\tilde K)(s - 1) - 1$ or a cycle $P_1$ such that $x_3^1, x_3^2 \in V(P_1)$,
$e(P_1) \le 2a_1$ and $v(P_1) = e(P_1)(s - 1)$. Since $Q(x_3^1)$ is true, there exists a $(2a_2 + a_3 + 1)$-cyclic extension $P_2$ of $(\{x_3^1\}, \varnothing)$ in $G$. Since $Q(x_3^2)$ is false, $P_2 \neq P_1$. Therefore, for the hypergraph $\tilde K = P_1 \cup P_2$ we have $e(\tilde K) \le 2a_1 + 2a_2 + a_3 + 1 = 2^{k - s + 1} + a$ and $v(\tilde K) \le e(\tilde K) (s - 1) - 1$.

So we have $\Pr[G^s(n, n^{-\alpha}) \models L] \le \Pr[G^s(n, n^{-\alpha}) \models L_{\mathcal K}]$, where $L_{\mathcal K}$ is the property of containing
a subhypergraph isomorphic to one of the hypergraphs from the set $\mathcal K$, where $\mathcal K$ is a finite set of $s$-hypergraphs with density at least $1/\alpha$. By theorem 13 there exists $\lim_{n \to \infty} \Pr[G^s(n, n^{-\alpha}) \models L_{\mathcal K}] \in (0, 1)$.
Hence  $\lim \sup_{n \to \infty} \Pr[G^s(n, n^{-\alpha}) \models L] < 1$.

Therefore, $G^s(n, n^{-\alpha})$ does not follow the zero-one $k$-law. \\

For vertices $x_1, \dots, x_k$ in a hypergraph $\mathcal G$ denote by
$\tilde N(x_1, \dots, x_k)$ the set of all vertices $y \in V(\mathcal G)$ such that
for all distinct $i_1, \ldots i_{s - 1} \in \{1, \dots, k\}$ we have $\{y, x_{i_1}, \dots, x_{i_{s - 1}}\} \in E(\mathcal G)$.

\textbf{Proof of theorem 10}.
Let $L$ be a first-order property which is expressed by the formula $\exists x_1 \ldots \exists x_{k - 11} \mbox{ }(Q_1(x_1, \ldots, x_{k - 11}) \land Q_2(x_1, \ldots, x_{k - 11}))$, where the formulae $Q_1(x_1, \ldots, x_{k - 11})$ and $Q_2(x_1, \ldots, x_{k - 11})$ are defined below.

First, let us introduce some notation. Fix vertices $x_1, \dots, x_{k - 11}, z$. Set $X = \tilde N(x_1, \dots, x_{k - 11})$. Consider the following predicates:
\begin{multline*}
T_z(x) = \exists v \mbox{ } ((v \in \tilde N(x_1, \ldots, x_{k - 13}, z, x)) \land \\ (\forall y \in X \mbox{ } ((y \neq x) \Rightarrow 
\bigcap_{\{a_1, \dots, a_{s - 2}\} \subset \{x_1, \ldots, x_{k - 13}, z\}}  (\neg N(v, y, a_1, \dots, a_{s - 2}))))),
\end{multline*}
\begin{multline*}
R_z(a, b) = \exists v \mbox{ } ((v \in\tilde N(x_1, \ldots, x_{k - 14}, z, a, b)) \land \\ (\forall y \in X
\mbox{ } (((y \neq a) \land (y \neq b)) \Rightarrow 
\bigcap_{\{a_1, \dots, a_{s - 2}\} \subset \{x_1, \ldots, x_{k - 14}, z\}}  (\neg N(v, y, a_1, \dots, a_{s - 2}))))).
\end{multline*}
Let $H_1(z) = \{x: T_z(x)\}$, $H_2(z) = \{x: \exists y \in H_1(z), R_z(x, y)\}$.
For a vertex $x$ set $N_z(x) = \{y: R_z(x, y)\}$.

Let 
$$
Q_1(x_1, \ldots, x_{k - 11}) = \exists z \mbox{ } (\varphi_1(z) \land \varphi_2(z) \land
\varphi_3(z) \land \varphi_4(z)),
$$
where
\begin{multline*}
\varphi_1(z) = \forall u_1 \in X \cap H_1(z) \mbox{ } \forall u_2 \in X \cap H_1(z) \mbox{ } ((u_1 \neq u_2) \Rightarrow ((N_z(u_1) \cap N_z(u_2) = \varnothing) \\ \land (N_z(u_1) \neq^* N_z(u_2)) \land (N_z(u_1) \cap H_1(z) = \varnothing))),
\end{multline*}

$$
\varphi_2(z) = \forall x \in X \mbox{ } ((x \in H_1(z)) \lor (x \in H_2(z))),
$$

$$
\varphi_3(z) = \forall u \in X \cap H_1(z) \mbox{ } (MIN_z(u) \Rightarrow (N_z(u) =^* H_1(z))),
$$

\begin{multline*}
\varphi_4(z) =  \forall u \in X \cap H_1(z) \mbox{ } \forall u_1 \in X \cap H_1(z) \mbox{ } \forall u_2 \in X \cap H_1(z) \\
((MIN_z(u) \land NEXT_z(u_1, u_2)) \Rightarrow (N_z(u_2) =^z N_z(u_1) \cdot^z N_z(u))),
\end{multline*}

$$
MIN_z(u) = \forall v \in X \cap H_1(z) \mbox{ } (N_z(v) \le^* N_z(u)),
$$

\begin{multline*}
NEXT_z(u_1, u_2) = (N_z(u_1) <^* N_z(u_2)) \land (\forall v \in X \cap H_1(z) \\ ((v \neq u_2) \land
(N_z(u_1) <^* N_z(v)) \Rightarrow (N_z(u_2) <^* N_z(v)))).
\end{multline*}

The formula $(A \le^* B)$ expresses the property that implies the inequality $|A| \le |B|$:
\begin{multline*}
(A \le^* B) = \exists \tilde z \mbox{ } (\forall x \mbox{ } ((x \in A \setminus B) \Rightarrow \exists y 
\mbox{ } ((y \in B \setminus A) \land R_{\tilde z}(x, y))) \land \\ 
((x \in B \setminus A) \Rightarrow (\neg (\exists y \mbox{ } \exists \tilde y \mbox{ }
((y \neq \tilde y) \land R_{\tilde z}(x, y) \land R_{\tilde z}(x, \tilde y)))))).
\end{multline*}
Denote $(A =^* B) = (A \le^* B) \land (B \le^* A)$, $(A <^* B) = (A \le^* B) \land (\neg (A =^* B))$.

The formula $(A =^z  B \cdot^z C)$ expresses the property that implies the inequality $|A| = |B||C|$:
\begin{multline*}
(A =^z  B \cdot^z C) = (\forall a \in A \mbox{ } \exists b \in B \mbox{ } \exists c \in C \mbox{ }
(R_z(a, b) \land R_z(a, c) \land \\
(\forall \tilde b \in B \mbox{ } ((b \neq \tilde b) \Rightarrow (\neg R_z(a, \tilde b)))) \land
(\forall \tilde c \in C \mbox{ } ((c \neq \tilde c) \Rightarrow (\neg R_z(a, \tilde c))))))
\land \\
(\forall b \in B \mbox{ } \forall c \in C \mbox{ } \exists a \in A \mbox{ }
(R_z(a, b) \land R_z(a, c) \land (\forall \tilde a \in A \mbox{ } ((a \neq \tilde a) \Rightarrow 
((\neg R_z(a, \tilde b)) \lor (\neg R_z(a, \tilde c))))))).
\end{multline*}

Let 
$$
Q_2(x_1, \ldots, x_{k - 11}) = \forall y \mbox{ } (\tilde N(x_1, \ldots, x_{k - 11}) \ge^*
\tilde N(x_1, \ldots, x_{k - 12}, y)).
$$

Let us say that vertices $x_1, \dots, x_{k - 11}$ are \textit{moderate}, if
$$
|\tilde N(x_1, \ldots, x_{k - 11})| = \max_y |\tilde N(x_1, \ldots, x_{k - 12}, y)|.
$$

Set
$$
\alpha = \frac{1}{{{k - 11} \choose {s - 1}}} + \frac{k - 10}{{{k - 11} \choose {s - 1}} \Sigma},
$$
where $\Sigma = 4 \frac{m^{m + 1} - m}{m - 1}$, $m = j (k - 10)$, $j \in \mathbb N$.

Let $\tilde \Omega_n$ be the set of all hypergraphs $\mathcal G$ from $\Omega_n$ which satisfy the following properties.
\begin{itemize}
\item[1)] For any vertices $x_1, \dots, x_{k - 11}$ we have $|\tilde N(x_1, \dots, x_{k - 11})| \le \frac{k - 11}{k - 10} \Sigma$.
\item[2)]  For any vertices $x_1, \dots, x_{k - 12}$ there exists a vertex $x_{k - 11}$ such that $|\tilde N(x_1, \dots, x_{k - 11})| \ge \frac{\Sigma}{k - 10} - 1$.
\item[3)] There are no subhypergraphs with at most $m^{m + 1}$ vertices and density greater than $1/\alpha$.
\item[4)] For any strictly balanced pair $(G, H)$ such that $\rho(G, H) < \frac{1}{\alpha}$ and $v(G) \le m^{m + 1}$, any collection of $v(H)$ vertices has a $(G, H)$-extension in $\mathcal G$.
\end{itemize}

Theorems 12 and 14 imply that $\Pr[G^s(n, n^{-\alpha}) \in \tilde \Omega_n] \to 1$, $n \to \infty$.

Let a hypergraph $\mathcal G \in \tilde \Omega_n$ satisfy $L$.
Then there exist vertices $x_1, \dots, x_{k - 11}$ satisfying $Q_1(x_1, \dots, x_{k - 11}) \land Q_2(x_1, \dots, x_{k - 11})$. Then $x_1, \dots, x_{k - 11}$ are moderate and $|\tilde N(x_1, \dots, x_{k - 11})| = \frac{t^{t + 1} - t}{t - 1} + t$ for some $t \in \mathbb N$.
If $t > m$, then $|\tilde N(x_1, \dots, x_{k - 11})| \ge \frac{(m + 1)^{m + 2} - m - 1}{m} + m + 1 > \frac{k - 11}{k - 10} \Sigma$ which contradicts property 1) from the definition of $\tilde \Omega_n$.
If $t < m$, then $|\tilde N(x_1, \dots, x_{k - 11})| \le \frac{(m - 1)^{m} - m + 1}{m} + m - 1 < \frac{\Sigma}{k - 10} - 1$ which contradicts property 2) from the definition of $\tilde \Omega_n$.
Thus $t = m$. Note that $m$ predicates $T_z$ are necessary for constructing the set $H_1(z)$,
$\frac{m^{m + 1} - m}{m - 1}$ predicates $R_z$  are necessary for constructing the set $H_2(z)$ and $2 \frac{m^{m + 1} - m^2}{m - 1}$ predicates $R_z$  are necessary for representing all predicates $A =^z B \cdot^z C$. Therefore, $\mathcal G$ contains a copy of a hypergraph $G(x_1, \dots, x_{k - 11})$ with $k - 10 + \Sigma$ vertices and ${{k - 11} \choose {s - 1}} \Sigma$ edges.
The hypergraph $G(x_1, \dots, x_{k - 11})$ is strictly balanced and its density equals $\frac{1}{\alpha}$.

Let us show that if a hypergraph $\mathcal G \in \tilde \Omega_n$ contains $G(x_1, \dots, x_{k - 11})$ for some $x_1, \dots, x_{k - 11}$ then it satisfies $L$. Suppose that there exist
$x_1, \dots, x_{k - 11}$ such that $\mathcal G$ contains $G(x_1, \dots, x_{k - 11})$.
Then property 3) from the definition of $\tilde \Omega_n$ implies that $\max_{y \neq x_{k - 11}} |\tilde N(x_1, \dots, x_{k - 12}, y)| \le \frac{\Sigma}{k - 10}$.
Therefore, $\forall y \mbox{ } |\tilde N(x_1, \dots, x_{k - 11})| \ge |\tilde N(x_1, \dots, x_{k - 12}, y)|$. The property 4) from the definition of $\tilde \Omega_n$ implies $\tilde L(x_1, \dots, x_{k - 11})$ because we need $\max_{y \neq x_{k - 11}} |\tilde N(x_1, \dots, x_{k - 12}, y)| \le \frac{\Sigma}{k - 10}$ predicates $R_{\tilde z}$ for this. The property 4) also implies $L(x_1, \dots, x_{k - 11})$  because we need at most $m^{m - 1} $predicates $R_{\tilde z}$ to represent predicates $\le^*$. Hence $\mathcal G$ satisfies $L$.

By theorem 13 the probability that $G^s(n, n^{-\alpha})$ contains $G(x_1, \dots, x_{k - 11})$ for some $x_1, \dots, x_{k - 11}$ tends to a constant $c \in (0, 1)$. Therefore,  
$\lim_{n \to \infty} \Pr[G^s(n, n^{-\alpha}) \models L] = c$.

Letting $j \to \infty$ we obtain that $\frac{1}{{{k - 11} \choose {s - 1}}} \in (S_k)'$. \\

\textbf{Proof of theorem 11}. Let $l = l(k)$. From the conditions of the theorem it follows that 
$l \ge 2$ if $s = 2$ and $l \ge s - 1$ if $s \ge 3$, so $\frac{1}{{l \choose {s - 1}}} \in (0, s - 1)$.
 Furthermore,  $l < k - 2$. Set $t = k - l - 2$, then $1 \le t < l$. Let $m \in \mathbb N$,
$\alpha = \frac{l + m}{(l - t + m) {l \choose {s - 1}}}$. 

Let $L$ be a first-order property which is expressed by the formula $\exists x_1 \ldots \exists x_l \mbox{ } Q(x_1, \ldots, x_l)$ with quantifier depth $k$, where
$$
Q(x_1, \ldots, x_l) = Q_1(x_1, \ldots, x_l)  \land (\neg (\exists z_1 \ldots \exists z_t \mbox{ } (Q_2(x_1, \ldots, x_l, z_1, \ldots, z_t) \land Q_3(x_1, \ldots, x_l, z_1, \ldots, z_t)))),
$$
\begin{multline*}
 Q_1(x_1, \ldots, x_l)  = \exists y_1 \ldots \exists y_{k - l} \mbox{ } ((y_1 \in \tilde N(x_1, \dots, x_l)) \land \ldots \land (y_{k - l} \in \tilde N(x_1, \dots, x_l)) \\ \land Q_4(x_1, \ldots, x_l, y_1, \ldots, y_{k - l}))
\end{multline*}
is a formula describing a strictly balanced $s$-hypergraph on vertices $x_1, \dots, x_l, y_1, \dots, y_{k - l}$ which has ${l \choose {s - 1}} (l + 2)$ edges and contains edges of the form
$\{y_i, x_{i_1}, \dots, x_{i_{s - 1}}\}$, where $i \in \{1, \dots, k - l\}$ and $i_1, \dots, i_{s - 1} \in \{1, \dots, l\}$ are distinct (such an $s$-hypergraph exists, since ${l \choose {s - 1}} (l + 2) \le {k \choose s}$), and $Q_4(x_1, \ldots, x_l, y_1, \ldots, y_{k - l})$ is a formula describing the rest ${l \choose {s - 1}} (l + 2) - {l \choose {s - 1}} (k - l)$ edges of this $s$-hypergraph,

\begin{multline*}
 Q_2(x_1, \ldots, x_l, z_1, \ldots, z_t)  = \forall y \mbox{ } ((y \in \tilde N(x_1, \dots, x_l)) \Rightarrow (\exists v \mbox{ } ((v \in \tilde N(z_1, \ldots, z_t, y, x_{t + 2}, \ldots, x_l)  \land \\
\bigcap_{\{a_1, \dots, a_{s - 1}\} \subset \{z_1, \ldots, z_t, x_1, \ldots, x_l\}, 
\{a_1, \ldots, a_{s - 1}\} \cap \{x_1, \ldots, x_{t + 1}\} \neq \varnothing}  (\neg N(v, a_1, \dots, a_{s - 1}))
)))),
\end{multline*}
\begin{multline*}
 Q_3(x_1, \ldots, x_l,  z_1, \ldots, z_t)  = \bigcap_{i = t + 1}^{l} (\exists u \mbox{ } ((u \in \tilde N(z_1, \ldots, z_t, x_1, \ldots, x_{i - 1}, x_{i + 1}, \ldots, x_l)) \land \\
\bigcap_{\{a_1, \ldots, a_{s - 2}\} \subset \{z_1, \ldots, z_t, x_1, \ldots, x_{i - 1}, x_{i + 1}, \ldots x_l\}}  (\neg N(u, x_i, a_1, \dots, a_{s - 2})).
))
\end{multline*}

Note that if $t = l - 1$ then
\begin{multline*}
 Q_2(x_1, \ldots, x_l, z_1, \ldots, z_t)  = \forall y \mbox{ } ((y \in \tilde N(x_1, \dots, x_l)) \Rightarrow (\exists v \mbox{ } ((v \in \tilde N(z_1, \ldots, z_t, y)) \land \\
\bigcap_{\{a_1, \dots, a_{s - 1}\} \subset \{z_1, \ldots, z_t, x_1, \ldots, x_l\}, 
\{a_1, \ldots, a_{s - 1}\} \cap \{x_1, \ldots, x_l\} \neq \varnothing}  (\neg N(v, a_1, \dots, a_{s - 1}))
)).
\end{multline*}

Let $\tilde \Omega_n$ be the set of all hypergraphs $\mathcal G$ from $\Omega_n$ which satisfy the following properties.
\begin{itemize}
\item[1)] For any strictly balanced pair $(G, H)$ such that $\rho(G, H) < \frac{1}{\alpha}$ and $v(G) \le 2(l + m + 1)$,  any collection of $v(H)$ vertices has a $(G, H)$-extension in $\mathcal G$.
\item[2)] For any hypergraph $G$ with $\rho^{\max}(G) > \frac{1}{\alpha}$ and $v(G) \le 2(l + m + 1)$, there is no copy of $G$ in $\mathcal G$. 
\end{itemize}

Theorems 12 and 14 imply that $\Pr[G^s(n, n^{-\alpha}) \in \tilde \Omega_n] \to 1$, $n \to \infty$.

Suppose that $\mathcal G \in \tilde \Omega_n$. Let us show that for any vertices $x_1, \dots, x_l$ satisfying $Q(x_1, \dots, x_l)$ we have $|\tilde N(x_1, \dots, x_l)| = m$.
Let $\tilde N(x_1, \dots, x_l) = \{y_1, \dots, y_{\chi}\}$. Let us prove that $\chi \le m$. Suppose that $\chi > m$. Since $Q_1(x_1, \ldots, x_l)$ is true, there exist $k - l$ vertices among $y_1, \dots, y_{\chi}$ (w.l.o.g. $y_1, \ldots, y_{k - l}$) such that $\mathcal G|_{\{x_1, \dots, x_l, y_1, \dots, y_{k - l}\}}$ contains at least ${l \choose {s - 1}} (l + 2) = {l \choose {s - 1}} (k - l) + f$ edges, where $f = (l - t) {l \choose {s - 1}}$. Then the density of the subhypergraph $\mathcal G|_{\{x_1, \dots, x_l, y_1, \dots, y_{m + 1}\}}$ is at least
$\frac{{l \choose {s - 1}} (m + 1) + f}{l + m + 1} > \frac{1}{\alpha}$. This contradicts property 1) from the definition of $\tilde \Omega_n$.
Let us prove that $\chi \ge m$. Suppose that $\chi < m$.
By the definition of $\tilde \Omega_n$, in $\mathcal G$ there are distinct vertices $z_1, \dots, z_t, v_1, \dots, v_{\chi}, u_1, \dots, u_{l - t}$ such that 
for any $i \in \{1, \dots, \chi\}$ we have $v_i \in \tilde N(z_1, \ldots, z_t, y_i, x_{t + 2}, \ldots, x_l)$ and for any $j \in \{1, \dots, l - t\}$ we have 
$u_j \in N(z_1, \ldots, z_t, x_1, \ldots, x_{t +  j - 1}, x_{t + j + 1}, \ldots, x_l)$.
Indeed, in this case the pair $(\mathcal G|_{W}, \mathcal G|_U)$, where $U = \{x_1, \dots, x_l, y_1, \dots, y_{\chi}\}$, $W = U \cup \{z_1, \dots, z_t, v_1, \dots, v_{\chi}, u_1, \dots, u_{l - t}\}$, is strictly balanced with the density
$\frac{(l - t + \chi) {l \choose {s - 1}}}{l + \chi} < \frac{(l - t + m) {l \choose {s - 1}}}{l + m} = \frac{1}{\alpha}$. This contradicts $Q(x_1, \ldots, x_l)$. Thus $\chi = m$.

Let $\tilde H$ be an $s$-hypergraph  with $V(H) = \{x_1, \dots, x_l, y_1, \dots, y_{k - l}\}$ and ${l \choose {s - 1}} (l + 2)$ edges which was determined above in defining of the formula $Q_1(x_1, \dots, x_l)$.
Let $H$ be an $s$-hypergraph with a set of vertices $V(H) = \{x_1, \dots, x_l, y_1, \dots, y_m\}$ and a set of edges 
$$E(H) = \bigcup_{i = 1}^{m} \{\{y_i, x_{i_1}, \ldots, x_{i_{s - 1}}\}, \{i_1, \ldots, i_{s - 1}\} \subset \{1, \ldots, l\}\} \cup E(\tilde H),$$
$G \supset H$ be an $s$-hypergraph with a set of vertices 
\newline
$V(G) = V(H) \cup \{z_1, \dots, z_t, v_1, \dots, v_m, u_1, \dots, u_{l - t}\}$ and a set of edges \begin{multline*}
E(G) = E(H) \cup \bigcup_{i = 1}^{m} \{\{v_i, a_1, \dots, a_{s - 1}\}, \{a_1, \dots, a_{s - 1}\} \subset \{y_i, z_1, \dots, z_t, x_{t + 2}, \dots, x_l\}\} \cup \\
\bigcup_{i = 1}^{l - t} \{\{u_i, a_1, \dots, a_{s - 1}\}, \{a_1, \dots, a_{s - 1}\} \subset \{z_1, \dots, z_t, x_1, \dots, x_{t + i - 1}, x_{t + i + 1}, \dots, x_l\}\}.
\end{multline*}  
Then $H$ is strictly balanced, the pair $(G, H)$ is strictly balanced and $\rho(H) = \rho(G, H) = \frac{1}{\alpha}$. Let $\tilde L$ be the property that in $\mathcal G$ there exists a copy of $H$ such that no copy of $G$ contains it.

Let us prove that a hypergraph $\mathcal G \in \tilde \Omega_n$ satisfies $L$ if and only if it satisfies $\tilde L$. Suppose that $\mathcal G$ satisfies $L$. Then there exist vertices $x_1, \dots, x_l$ satisfying $Q(x_1, \dots, x_l)$. As we have shown, $|\tilde N(x_1, \dots, x_l)| = m$. Then $\mathcal G|_{\{x_1, \dots, x_l\} \cup \tilde N(x_1, \dots, x_l)}$ is isomorphic to $H$. The existence of a copy of $G$ containing $\mathcal G|_{\{x_1, \dots, x_l\} \cup \tilde N(x_1, \dots, x_l)}$ contradicts $Q(x_1, \dots, x_l)$. Thus $\mathcal G$ satisfies $\tilde L$.

Suppose that $\mathcal G$ satisfies $\tilde L$. Then there exists a copy of $H$ such that no copy of $G$ contains it. Let $x_1, \dots, x_l, y_1, \dots, y_m$ be the vertices of this copy of $H$. Then $\tilde N(x_1, \dots, x_l) = \{y_1, \dots, y_m\}$. Indeed, otherwise $\mathcal G$ contains a subhypergraph with $m + l + 1$ vertices and density more than $1/\alpha$ which contradicts the definition of $\tilde \Omega_n$. Suppose that $Q(x_1, \dots, x_l)$ is false. Since $Q_1(x_1, \dots, x_l)$ is true, there exist vertices $z_1, \dots, z_t \in V(\mathcal G)$ such that the predicates $Q_2(x_1, \ldots, x_l, z_1, \ldots, z_t)$ and $Q_3(x_1, \ldots, x_l, z_1, \ldots, z_t)$ are true. The predicates $Q_2(x_1, \ldots, x_l, z_1, \ldots, z_t)$ and $Q_3(x_1, \ldots, x_l, z_1, \ldots, z_t)$ imply the existence of a number $r \le m$, distinct vertices $v_1, \dots, v_r$, $u_1, \dots, u_{l - t}$ and disjoint subsets $Y_1, \dots, Y_r$, $Y_1 \sqcup \ldots \sqcup Y_r = \{y_1, \dots, y_m\}$, such that 
\begin{itemize}
\item[1)] for any $i \in \{1, \dots, r\}$ and for any $y \in Y_i$ we have $v_i \in \tilde N(z_1, \ldots, z_t, y, x_{t + 2}, \ldots, x_l)$,
\item[2)] for any $j \in \{1, \dots, l - t\}$ we have $u_j \in \tilde N(z_1, \ldots, z_t, x_1, \ldots, x_{t +  j - 1}, x_{t + j + 1}, \ldots, x_l)$.
\end{itemize}
Set $W = \{x_1, \dots, x_l, y_1, \dots, y_m\} \cup \{z_1, \dots, z_t, v_1, \dots, v_r, u_1, \dots, u_{l - m}\}$. Then the density of $\mathcal G|_W$ is at least 
$$\frac{2(l - t + m) {l \choose {s - 1}} - (m - r) {{l - 1} \choose {s - 1}}}{2l + m + r},$$
which is greater than $1/\alpha$  if $r < m$. The definition of $\tilde \Omega_n$ implies that $r = m$ and $\mathcal G|_W$ is a copy of $G$. This contradicts the property $\tilde L$. Therefore, $Q(x_1, \dots, x_l)$ is true, and hence $\mathcal G$ satisfies $L$.

By proposition 1, there exists $\lim_{n \to \infty} \Pr[G^s(n, n^{-\alpha})\models \tilde L] = c \in (0, 1)$. Since 
\newline
$\lim_{n \to \infty} \Pr[G^s(n, n^{-\alpha}) \in \tilde \Omega_{n}] = 1$, it follows from the above that 
$\lim_{n \to \infty} \Pr[G^s(n, n^{-\alpha}) \models L] = c.$ Letting $m \to \infty$ we obtain that $\frac{1}{{l \choose {s - 1}}} \in (S_k)'$.

\section{Proofs of propositions}
\indent

\textbf{Proof of proposition 2}. 

1) Let $G_0, G_1, \ldots, G_t, G$ be the $m$-decomposition of $G$. Set 
$e_i = e(G_i, G_{i - 1}) - 1$, $v_i = v(G_i, G_{i - 1}) - e_i (s - 1)$ for all $i \in \{1, \ldots, t\}$. 
Set $e_0 = e(G) - e(G_t)$. Note that $e_i \le m - 1$, $0 \le v_i \le s - 2$.
We have
$$
\frac{1}{\rho(G)} = \frac{\sum_{i = 1}^{t} (e_i (s - 1) + v_i) + 1}{\sum_{i = 1}^{t} e_i + t + e_0} = s - 1 - \frac{(t + e_0)(s - 1) - \sum v_i - 1}{\sum_{i = 1}^{t} e_i + t + e_0}.
$$
Note that $(t + e_0)(s - 1) - \sum v_i - 1 \ge (t + e_0)(s - 1) - t(s - 2) - 1 = t + e_0(s - 1) - 1$. Therefore, $\frac{1}{\rho(G)} \le s - 1$ and
the equality holds if and only if $e_0 = 0$, $t = 1$ and $v_1 = s - 2$. Hence either $G$ is a cyclic $m$-extension of $(\{x\}, \varnothing)$ and $\frac{1}{\rho^{\max}(G)} = s - 1$
or $\frac{1}{\rho(G)} < s - 1$. 

Suppose that $\frac{1}{\rho(G)} < s - 1$. Then 
$$\frac{1}{\rho(G)} = s - 1 - \frac{1}{\tau}, \quad \mbox{where} \quad
 \tau =\frac{\sum_{i = 1}^{t} e_i + t + e_0}{(t + e_0)(s - 1) - \sum v_i - 1}.$$
Since $\rho(G) < \frac{m}{m(s - 1) - 1}$, $\tau \ge m$. We have $\tau = m + \frac{a}{b}$, where
$$a = \sum_{i = 1}^{t} e_i + t + e_0 - m(s - 1)(e_0 + t) + m\left(\sum v_i + 1\right), \quad b = (t + e_0)(s - 1) - \sum v_i - 1.$$
Since $e_i \le m - 1$ and $v_i \le s - 2$, we obtain
$$
a \le mt + e_0 - m(s - 1)(e_0 + t) + m ((s - 2)t + 1) \le m - e_0 (m(s - 1) - 1) \le m.
$$ 

Let $H \subset G$, $\rho(G) < \rho(H) < \frac{m}{m(s - 1) - 1}$. Let $v = v(G) - v(H)$, $e = e(G) - e(H)$. Then
$$
\frac{1}{\rho(H)} = \frac{\sum_{i = 1}^{t} (e_i (s - 1) + v_i) + 1 - v}{\sum_{i = 1}^{t} e_i + t + e_0 - e} = s - 1 - \frac{1}{\tau'},
$$
where
$$
\tau' = \frac{\sum_{i = 1}^{t} e_i + t + e_0 - e}{(t + e_0)(s - 1) - \sum v_i - 1 + v - e(s - 1)}.
$$
Since $\rho(H) < \frac{m}{m(s - 1) - 1}$, $\tau' \ge m$. We have $\tau' = m + \frac{a'}{b'}$, where
$$a' = a + e (m(s - 1) - 1) - m v, \quad b' = b + v - e(s - 1).$$
Since $G$ is connected, $v \le e (s - 1)$. Therefore,  $b' \le b$. Since $\rho(G) < \rho(H)$, we have $\frac{a'}{b'} \le \frac{a}{b}$ and $a' < a \le m$.
So we obtain that $\frac{1}{\rho(H)} = s - 1 - \frac{1}{m + a'/b'}$, where $a' \le m$.

2) Let $G_0, G_1, \ldots, G_t, G$ be the $m$-decomposition of a hypergraph $G \in \mathcal H_m$. 
Then $\frac{1}{\rho(G)} = s - 1 - \frac{1}{\tau}$, where
$$
\tau \le \frac{m t}{t(s - 1) - t(s - 2) - 1}= \frac{m t}{t - 1}.
$$
Note that $|V(G)| \le t (m - 1) (s - 1) + t (s - 2) + 1$. Therefore, if $|V(G)| \ge l (m - 1) (s - 1) + l (s - 2) + 1$, then 
$\frac{1}{\rho(G)} \le s - 1 - \frac{l - 1}{ml}$. For any $\rho < \frac{m}{m(s - 1) - 1}$, there exists $l = l(\rho)$ such that
$s - 1 - \frac{l - 1}{ml} < 1/\rho$.  Then $\rho(G) > \rho$ for any $G \in \mathcal H_m$ with $|V(G)| \ge  l(\rho) ((m - 1) (s - 1) + (s - 2)) + 1$. Since $v(G_{i + 1}) - v(G_i) \le (m - 1)(s - 1) + s - 2$,  any hypergraph $H \in \mathcal H_m$ with at least
$ (l(\rho) + 1) ((m - 1) (s - 1) + (s - 2)) + 1$ vertices contains a subhypergraph $H' \in \mathcal H_m$ with $v(H') \in 
\{l(\rho) ((m - 1) (s - 1) + (s - 2)) + 1, \ldots, (l(\rho) + 1) ((m - 1) (s - 1) + (s - 2)) + 1\}$.
Therefore, $\eta = \eta(\rho) =  (l(\rho) + 1) ((m - 1) (s - 1) + (s - 2)) + 1$ satisfies the statement of proposition. \\

\textbf{Proof of proposition 3}.
Consider a pair $(\tilde G \cup \tilde H, \tilde G)$, where $V(\tilde H) \cap V(\tilde G) = \{x\}$,
$V(\tilde H) = \{x, y_1, \ldots, y_{m(s - 1)}\}$, 
$$E(\tilde H) = \{\{x, y_1, \ldots, y_{s - 1}\}\} \cup \{\{y_{i(s - 1)}, \ldots, y_{(i + 1)(s - 1)}\}: 1 \le i < m\}.$$
Then $\rho(\tilde G \cup \tilde H, \tilde G) = \frac{1}{s - 1}$ and $(\tilde G \cup \tilde H, \tilde G)$ is $1/\rho$-safe.
Since $G$ is $(n_1(\rho), n_2(\rho), n_3, n_4)$-sparse, there exists a subhypergraph $G_1 \subset G$ such that $G_1$ is an exact $(\tilde G \cup \tilde H, \tilde G)$-extension of $\tilde G$
and $(G_1, \tilde G)$ is $(K_1, K_2)$-maximal in $G$ for all $(K_1, K_2) \in \mathcal K^{\rho}$. Then there exists a vertice $y \in V(G_1)$ such that $d_{G_1}(y, \tilde G) = m$. Moreover, $G_1 \setminus (\tilde G \setminus \{x\})$ is a path which connects $x$ and $y$ 
and has $m$ edges and $m(s - 1) + 1$ vertices. Since $(G_1, \tilde G)$ is $(K_1, K_2)$-maximal in $G$ for all $(K_1, K_2) \in \mathcal K^{\rho}$, there does not exist a path $P' \neq P$ which connects $y$ with a vertex from $\tilde G$ and has at most $m$ edges. \\

\textbf{Proof of proposition 4}.
Let $X_i$ and $Y_i$ be the hypergraphs chosen in the $i$-th round by Spoiler and Duplicator respectively. We denote vertices which are chosen in the first $i$ rounds in $X_i$ and $Y_i$ by $x^1_i, \ldots, x^i_i$ and $y^1_i, \ldots, y^i_i$.

Let $i$ rounds be finished, where $1 \le i \le k - s + 2$. Let $r \in \{1, \ldots, i\}$. Let $\tilde X^1_i, \ldots, \tilde X^r_i \subset X_i$ and $\tilde Y^1_i, \ldots, \tilde Y^r_i \subset Y_i$ be subhypergraphs of $X_i$ and $Y_i$ respectively. We say that $\tilde X^1_i, \ldots, \tilde X^r_i$ and $\tilde Y^1_i, \ldots, \tilde Y^r_i$ are \textit{$(k, i, r)$-regular equivalent} in $(X_i, Y_i)$, if the following properties hold.
\begin{itemize}
\item[(I)] $x^1_i, \ldots, x^i_i \in V(\tilde X^1_i \cup \ldots \cup \tilde X^r_i)$, $y^1_i, \ldots, y^i_i \in V(\tilde Y^1_i \cup \ldots \cup \tilde Y^r_i).$
\item[(II)] For any distinct $j_1, j_2 \in \{1, \ldots, r\}$, the inequalities $d_{X_i}(\tilde X^{j_1}_i, \tilde X^{j_2}_i) > 2^{k - i - s + 2}$,
$d_{Y_i}(\tilde Y^{j_1}_i, \tilde Y^{j_2}_i) > 2^{k - i - s + 2}$ hold.
\item[(III)] For any $j \in \{1, \ldots, r\}$, there is no cyclic $2^{k - i - s + 2}$-extension of $\tilde X^j_i$ in the hypergraph $X_i$
and there is no cyclic $2^{k - i - s + 2}$-extension of $\tilde Y^j_i$ in the hypergraph $Y_i$.
\item[(IV)] Cardinalities of the sets $V(\tilde X^1_i \cup \ldots \cup \tilde X^r_i)$ and $V(\tilde Y^1_i \cup \ldots \cup \tilde Y^r_i)$ are at most $\eta(\rho) + (i - 1) 2^{k - s + 1}$.
\item[(V)] The hypergraphs $\tilde X^j_i$ and $\tilde Y^j_i$ are isomorphic for any $j \in \{1, \ldots, r\}$ and there exists a corresponding isomorphism (one for all these pairs of hypergraphs) which maps the vertices $x^l_i$ to the vertices $y^l_i$ for all $l \in \{1, \ldots, i\}$.
\end{itemize}

The main idea of Duplicator's strategy is the following. Duplicator should play in such way that after $i$ rounds, where $1 \le i \le k - s + 2$, there exist a number $r \in \{1, \ldots, i\}$ and $(k, i, r)$-regular equivalent hypergraphs 
$\tilde X^1_i, \ldots, \tilde X^r_i \subset X_i$ and $\tilde Y^1_i, \ldots, \tilde Y^r_i \subset Y_i$.

Let us describe Duplicator's strategy in the first round. By virtue of lemma 1 and $(n_1(\rho), n_2(\rho),$ $n_3, n_4)$-sparseness of $X_1$ it does not contain any hypergraph $H \in \mathcal H_{2^{k - s + 1}}$ with $v(H) \ge \eta(\rho)$. Let $\tilde X^1_1$ be the subhypergraph of $X_1$ with maximum number of vertices such that $x^1_1 \in \tilde X^1_1$ and $\tilde X^1_1 \in \mathcal H_{2^{k - s + 1}}$. Then $v(\tilde X^1_1) < \eta(\rho) < n_1(\rho)$. Since $X_1$ is $(n_1(\rho), n_2(\rho), n_3, n_4)$-sparse, the density $\rho(\tilde X^1_1) \le \rho$. Proposition 2 implies that $\rho(\tilde X^1_1) \neq \rho$, since 
$\rho \notin \mathcal Q_k$. Therefore, $\rho(\tilde X^1_1) < \rho$ and by the property of $(n_1(\rho), n_2(\rho), n_3, n_4)$-sparseness $Y_1$ contains a induced subhypergraph $\tilde Y^1_1$ which is isomorphic to $\tilde X^1_1$ and $(K_1, K_2)$-maximal in $Y_1$ for all $(K_1, K_2) \in \mathcal K^{\rho}$. Let $\varphi_1 \colon \tilde X^1_1 \to \tilde Y^1_1$ be an isomorphism. Then Duplicator chooses the vertex $y^1_1 = \varphi_1(x^1_1)$. Since $\tilde Y^1_1$ is
$(K_1, K_2)$-maximal in $Y_1$ for all $(K_1, K_2) \in \mathcal K^{\rho}$, there is no cyclic $2^{k - s + 1}$-extension of $\tilde Y^1_1$ in $Y_1$.
Therefore, $\tilde X^1_1$ and $\tilde Y^1_1$ are $(k, 1, 1)$-regular equivalent.

Let $i$ rounds be finished, where $1 \le i \le k - s + 1$. Let us describe Duplicator's strategy in the $(i + 1)$-th round.
If $X_{i + 1} = X_i$, then set $\tilde X^j_{i + 1} = \tilde X^j_i$, $\tilde Y^j_{i + 1} = \tilde Y^j_i$.
Otherwise, set $\tilde X^j_{i + 1} = \tilde Y^j_i$, $\tilde Y^j_{i + 1} = \tilde X^j_i$. Let $\varphi_{i + 1} \colon \tilde X_{i + 1}^1 \cup \ldots \cup \tilde X_{i + 1}^r \to \tilde Y_{i + 1}^1 \cup \ldots \cup \tilde Y_{i + 1}^r$ be an isomorphism such that $\varphi_{i + 1}(\tilde X^j_{i + 1}) = \tilde Y^j_{i + 1}$ for all $j \in \{1, \ldots, r\}$ and $\varphi_{i + 1}(x^l_{i + 1}) = y^l_{i + 1}$ for all $l \in \{1, \ldots, i\}$.

Consider three cases.

1) Spoiler chooses a vertex $x^{i + 1}_{i + 1} \in V(\tilde X^1_{i + 1} \cup \ldots \cup \tilde X^{r}_{i + 1})$. Then Duplicator chooses
$y^{i + 1}_{i + 1} = \varphi_{i + 1}(x^{i + 1}_{i + 1})$. Obviously, hypergraphs $\tilde X^1_{i + 1}, \ldots, \tilde X^r_{i + 1}$ and
$\tilde Y^1_{i + 1}, \ldots, \tilde Y^r_{i + 1}$ are $(k, i + 1, r)$-regular equivalent.

2) Spoiler chooses a vertex $x^{i + 1}_{i + 1} \notin V(\tilde X^1_{i + 1} \cup \ldots \cup \tilde X^{r}_{i + 1})$ such that
$d_{X_{i + 1}}(x^{i + 1}_{i + 1}, \tilde X^1_{i + 1} \cup \ldots \cup \tilde X^{r}_{i + 1}) \le 2^{k - i - s + 1}$. Then by (II) there exists exactly one $j \in \{1, \ldots, r\}$ such that $d_{X_{i + 1}}(x^{i + 1}_{i + 1}, \tilde X^j_{i + 1}) \le 2^{k - i - s + 1}$. By (III) there exists exactly one path $C_{X_{i + 1}}$ of length at most $2^{k - i - s + 1}$ connecting $x^{i + 1}_{i + 1}$ with a vertex from $\tilde X^j_{i + 1}$. Let $\tilde x^{j}_{i + 1}$ be the vertex of $C_{X_{i + 1}}$ belonging to $\tilde X^j_{i + 1}$.
(III) also implies that $v(C_{X_{i + 1}}) = e(C_{X_{i + 1}})(s - 1) + 1$ (otherwise $C_{X_{i + 1}}$ gives a cyclic $2^{k - i - s + 2}$-extension of $\tilde X^j_{i + 1}$). 
By (IV) $|V(\tilde Y^1_{i + 1} \cup \ldots \cup \tilde Y^r_{i + 1})| \le \eta(\rho) + (i - 1) 2^{k - s + 1} \le n_2(\rho)$.
Therefore, by Proposition 3 there exist a vertex $y^{i + 1}_{i + 1} \in Y_i$ and a path $C_{Y_{i + 1}} \subset Y_i$ 
such that $d_{Y_{i + 1}}(y^{i + 1}_{i + 1}, \tilde Y^1_{i + 1} \cup \ldots \cup \tilde Y^r_{i + 1}) = d_{X_{i + 1}}(x^{i + 1}_{i + 1}, \tilde X^1_{i + 1} \cup \ldots \cup \tilde X^r_{i + 1})$, 
$C_{Y_{i + 1}}$ is isomorphic to $C_{X_{i + 1}}$ and connects $y^{i + 1}_{i + 1}$ with $\varphi_{i + 1}(\tilde x^{j}_{i + 1})$ and there does not exist any path $C' \subset Y_{i + 1}$, $C' \neq C_{Y_{i + 1}}$, such that $C'$ connects $y^{i + 1}_{i + 1}$
with a vertex from $ \tilde Y^1_{i + 1} \cup \ldots \cup \tilde Y^r_{i + 1}$ and $e(C') \le e(C_{Y_{i + 1}})$.
Then there exists an isomorphism $C_{X_{i + 1}} \cup \tilde X^1_{i + 1} \cup \ldots \cup \tilde X^r_{i + 1} \to C_{Y_{i + 1}} \cup \tilde Y^1_{i + 1} \cup \ldots \cup \tilde Y^r_{i + 1}$ which maps the vertices $x^1_{i + 1}, \ldots, x^{i + 1}_{i + 1}$ to
$y^1_{i + 1}, \ldots, y^{i + 1}_{i + 1}$ respectively.
Redefine the hypergraphs $\tilde X^{j}_{i + 1}$, $\tilde Y^j_{i + 1}$: $\tilde X^j_{i + 1} = \tilde X^j_{i + 1} \cup C_{X_{i + 1}}$, $\tilde Y^j_{i + 1} = \tilde Y^j_{i + 1} \cup C_{Y_{i + 1}}$. 
Then $\tilde X^1_{i + 1}, \ldots, \tilde X^r_{i + 1}$ and $\tilde Y^1_{i + 1}, \ldots, \tilde Y^r_{i + 1}$ satisfy (I) and (V). Let us prove that properties (II), (III), (IV) are also satisfied.

Show that (II) holds. It is sufficient to show that $d_{X_{i + 1}}(\tilde X^{j}_{i + 1}, \tilde X^{j_2}_{i + 1}) > 2^{k - i - s + 1}$ and 
\newline
$d_{Y_{i + 1}}(\tilde Y^{j}_{i + 1}, \tilde Y^{j_2}_{i + 1}) > 2^{k - i - s + 1}$ for any $j_2 \neq j$.
Suppose that $d_{X_{i + 1}}(\tilde X^j_{i + 1}, \tilde X^{j_2}_{i + 1}) \le 2^{k - i - s + 1}$. 
Then there exists a vertex $u \in V(C_{X_{i + 1}})$ such that $d_{X_{i + 1}}(u, \tilde X^{j_2}_{i + 1}) \le 2^{k - i - s + 1}$.
Hence $d_{X_{i + 1}}(\tilde X^j_{i + 1} \setminus (C_{X_{i + 1}} \setminus \{\tilde x^j_{i + 1}\}), \tilde X^{j_2}_{i + 1}) \le 2^{k - i - s + 2}$. This contradicts the fact that $\tilde X^1_i, \ldots, \tilde X^r_i$ and $\tilde Y^1_i, \ldots, \tilde Y^r_i$ satisfy (II).
Therefore, $d_{X_{i + 1}}(\tilde X^j_{i + 1}, \tilde X^{j_2}_{i + 1}) > 2^{k - i - s + 1}$. The inequality $d_{Y_{i + 1}}(\tilde Y^j_{i + 1}, \tilde Y^{j_2}_{i + 1}) > 2^{k - i - s + 1}$ is proved analogously.

Show that (III) holds. It is sufficient to show that there is no cyclic $2^{k - i - s + 1}$-extension of $\tilde X^j_{i + 1}$ in $X_{i + 1}$ and there is no cyclic $2^{k - i - s + 1}$-extension of $\tilde Y^j_{i + 1}$ in $Y_{i + 1}$.
Suppose that there exists $W \subset X_{i + 1}$ which is a cyclic $2^{k - i - s + 1}$-extension of $\tilde X^j_{i + 1}$.
Then $e(W, \tilde X^j_{i + 1} \setminus (C_{X_{i + 1}} \setminus \{\tilde x^j_{i + 1}\})) \le 2^{k - i - s + 2}$ and there exists a subhypergraph $W' \subset W$ which is a cyclic $2^{k - i - s + 2}$-extension of $\tilde X^j_{i + 1} \setminus (C_{X_{i + 1}} \setminus \{\tilde x^j_{i + 1}\})$.
This contradicts the fact that $\tilde X^1_i, \ldots, \tilde X^r_i$ and $\tilde Y^1_i, \ldots, \tilde Y^r_i$ satisfy (III).
Therefore, there is no cyclic $2^{k - i - s + 1}$-extension of $\tilde X^j_{i + 1}$ in $X_{i + 1}$. The non-existence of cyclic $2^{k - i - s + 1}$-extensions of $\tilde Y^j_{i + 1}$ in $Y_{i + 1}$ is proved analogously.

Show that (IV) holds. 
We have
$$
|V(\tilde X^1_{i + 1} \cup \ldots \cup \tilde X^{r}_{i + 1})| \le |V(\tilde X^1_{i} \cup \ldots \cup \tilde X^{r}_{i})| + 2^{k - i - s + 1} \le \eta(\rho) + i 2^{k - s + 1},
$$
therefore, (IV) holds.

3) Spoiler chooses a vertex $x^{i + 1}_{i + 1} \in X_{i + 1}$ such that
$d_{X_{i + 1}}(x^{i + 1}_{i + 1}, \tilde X^1_{i + 1} \cup \ldots \cup \tilde X^{r}_{i + 1}) > 2^{k - i - s + 1}$.
 Let $\tilde X^{r + 1}_{i + 1}$ be the subhypergraph of $X_{i + 1}$ with maximum number of vertices such that $x^{i + 1}_{i + 1} \in \tilde X^{r + 1}_{i + 1}$ and $\tilde X^{r + 1}_{i + 1} \in \mathcal H_{2^{k - i - s + 1}}$.
By Proposition 2 either $\tilde X^{r + 1}_{i + 1}$ is a cyclic $2^{k - i - s + 1}$-extension of $(\{x^{i + 1}_{i + 1}\}, \varnothing)$ with $1/\rho(\tilde X^{r + 1}_{i + 1}) = s - 1$ or $1/\rho\tilde (X^{r + 1}_{i + 1}) = s - 1 - \frac{1}{2^{k - i - s + 1} + a/b}$,
 where $a, b \in \mathbb N$, $a \le 2^{k - i - s + 1}$.  Note that $s - 1 - \frac{1}{2^{k - i - s + 1} + a/b} \le s - 1 - \frac{1}{2^{k - i - s + 2}} \le s - 1 - \frac{1}{2^{k - s + 1}} < 1/\rho$. Since $X_{i + 1}$ is $(n_1(\rho), n_2(\rho), n_3, n_4)$-sparse,
we obtain that $\tilde X^{r + 1}_{i + 1}$ is a cyclic $2^{k - i - s + 1}$-extension of $\{x^{i + 1}_{i + 1}\}$ and $1/\rho(\tilde X^{r + 1}_{i + 1}) = s - 1$.
By (IV) $|V(\tilde Y^1_{i + 1} \cup \ldots \cup \tilde Y^r_{i + 1})| \le \eta(\rho) + (i - 1) 2^{k - s + 1} \le n_2(\rho)$.
Therefore, by Proposition 3 there exists a vertex $y^{i + 1}_{i + 1} \in Y_{i + 1}$ such that $d_{Y_{i + 1}}(y^{i + 1}_{i + 1}, \tilde Y^1_{i + 1} \cup \ldots \cup \tilde Y^r_{i + 1}) = 2^{k - i - s + 1} + 1$.
Since $Y_i$ is $(n_1(\rho), n_2(\rho), n_3, n_4)$-sparse, there exists a subhypergraph $\tilde Y^{r + 1}_{i + 1} \subset Y_{i + 1}$ such that $\tilde Y^{r + 1}_{i + 1}$ is an exact $(\tilde X^{r + 1}_{i + 1}, (\{x^{i + 1}_{i + 1}\}, \varnothing))$-extension
of $(\{y^{i + 1}_{i + 1}\}, \varnothing)$ and $\tilde Y^{r + 1}_{i + 1}$ is $(K_1, K_2)$-maximal in $Y_{i + 1}$ for all $(K_1, K_2) \in \mathcal K^{\rho}$.

Let us prove that  $\tilde X^1_{i + 1}, \ldots, \tilde X^{r + 1}_{i + 1}$ and $\tilde Y^1_{i + 1}, \ldots, \tilde Y^{r + 1}_{i + 1}$ satisfy (I)-(V). Obsiously, (I) and (V) hold.

Show that (II) holds. It is sufficient to show that $d_{X_{i + 1}}(\tilde X^{r + 1}_{i + 1}, \tilde X^j_{i + 1}) > 2^{k - i - s + 1}$ and 
\newline
$d_{Y_{i + 1}}(\tilde Y^{r + 1}_{i + 1}, \tilde Y^j_{i + 1}) > 2^{k - i - s + 1}$ for any $j \in \{1, \ldots, r\}$.
Suppose that $d_{X_{i + 1}}(\tilde X^{r + 1}_{i + 1}, \tilde X^j_{i + 1}) \le 2^{k - i - s + 1}$. Then there exists a vertex $u \in V(\tilde X^{r + 1}_{i + 1})$ such that $d_{X_{i + 1}}(u, \tilde X^j_{i + 1}) \le 2^{k - i - s + 1}$. Note that $\tilde X^{r + 1}_{i + 1} \cap \tilde X^j_{i + 1} = \varnothing$,
since $d_{X_{i + 1}}(x, \tilde X^1_{i + 1} \cup \ldots \cup \tilde X^r_{i + 1}) > 2^{k - i - s + 1}$.
Since $\tilde X^{r + 1}_{i + 1}$ is
a cyclic $2^{k - i - s + 1}$-extension of $(\{x^{i + 1}_{i + 1}\}, \varnothing)$, 
$\tilde X^{r + 1}_{i + 1} \cap \tilde X^j_{i + 1} = \varnothing$ and 
$d_{X_{i + 1}}(u, \tilde X^j_{i + 1}) \le 2^{k - i - s + 1}$, we obtain that there exists a cyclic $2^{k - i - s + 2}$-extension of $\tilde X^j_{i + 1}$ in $X_{i + 1}$. This is a contradiction.
Therefore, $d_{X_{i + 1}}(\tilde X^{r + 1}_{i + 1}, \tilde X^j_{i + 1}) > 2^{k - i - s + 1}$.  
The inequality $d_{Y_{i + 1}}(\tilde Y^{r + 1}_{i + 1}, \tilde Y^{j}_{i + 1}) > 2^{k - i - s + 1}$ is proved analogously.

Show that (III) holds. It is sufficient to show that there is no cyclic $2^{k - i - s + 1}$-extension of $\tilde X^{r + 1}_{i + 1}$ in $X_{i + 1}$ and 
 there is no cyclic $2^{k - i - s + 1}$-extension of $\tilde Y^{r + 1}_{i + 1}$ in $Y_{i + 1}$. By the definition of $\tilde X^{r + 1}_{i + 1}$ there is no $2^{k - i - s + 1}$-extension of $\tilde X^{r + 1}_{i + 1}$ in $X_{i + 1}$.
The hypergraph $\tilde Y^{r + 1}_{i + 1}$ does not contain any $2^{k - i - s + 1}$-extension of $\tilde Y^{r + 1}_{i + 1}$, since $\tilde Y^{r + 1}_{i + 1}$ is $(K_1, K_2)$-maximal in $Y_{i + 1}$ for all $(K_1, K_2) \in \mathcal K^{\rho}$.

Show that (IV) holds.
We have
$$
|V(\tilde X^1_{i + 1} \cup \ldots \cup \tilde X^{r + 1}_{i + 1})| \le |V(\tilde X^1_{i} \cup \ldots \cup \tilde X^{r}_{i})| + 2^{k - i - s + 1} \le \eta(\rho) + i 2^{k - s + 1},
$$
therefore, (IV) holds.

After $k - s + 2$ rounds we have vertices $ x^1_{k - s + 2}, \ldots, x^{k - s + 2}_{k - s + 2} \in X_{k - s + 2}$, $y^1_{k - s + 2}, \ldots, y^{k - s + 2}_{k - s + 2} \in Y_{k - s + 2}$ and hypergraphs
$\tilde X^1_{k - s + 2}, \ldots, \tilde X^r_{k - s + 2}$, $\tilde Y^1_{k - s + 2}, \ldots, \tilde Y^r_{k - s + 2}$ satisfying (I)-(V). 

Let us describe the strategy of Duplicator in the $(i + 1)$-th round, where $k - s + 2 \le i \le k - 1$.
If $X_{i + 1} = X_{k - s + 2}$, then set $\tilde X^j_{i + 1} = \tilde X^j_{k - s + 2}$, $\tilde Y^j_{i + 1} = \tilde Y^j_{k - s + 2}$.
Otherwise, set $\tilde X^j_{i + 1} = \tilde Y^j_{k - s + 2}$, $\tilde Y^j_{i + 1} = \tilde X^j_{k - s + 2}$. Let $\varphi_{i + 1} \colon \tilde X_{i + 1}^1 \cup \ldots \cup \tilde X_{i + 1}^r \to \tilde Y_{i + 1}^1 \cup \ldots \cup \tilde Y_{i + 1}^r$ be an isomorphism such that $\varphi_{i + 1}(\tilde X^j_{i + 1}) = \tilde Y^j_{i + 1}$ for all $j \in \{1, \ldots, r\}$ and $\varphi_{i + 1}(x^l_{i + 1}) = y^l_{i + 1}$ for all $l \in \{1, \ldots, k - s + 2\}$.
If Spoiler chooses a vertex $x^{i + 1}_{i + 1} \in V(\tilde X^1_{i + 1} \cup \ldots \cup \tilde X^{r}_{i + 1})$, then Duplicator chooses
$y^{i + 1}_{i + 1} = \varphi_{i + 1}(x^{i + 1}_{i + 1})$. Otherwise, Duplicator chooses an arbitrary vertex $y^{i + 1}_{i + 1} \notin V(\tilde Y^1_{i + 1} \cup \ldots \cup \tilde Y^{r}_{i + 1})$.

Let us show that Duplicator wins, i.e., the hypergraphs $X_{k}|_{\{x^1_k, \ldots, x^k_k\}}$ and $Y_{k}|_{\{y^1_k, \ldots, y^k_k\}}$ are isomorphic. Since $x^1_{k}, \ldots, x^{k - s + 2}_k \in V(\tilde X^1_{k} \cup \ldots \cup \tilde X^r_{k})$
and $y^1_{k}, \ldots, y^{k - s + 2}_k \in V(\tilde Y^1_{k} \cup \ldots \cup \tilde Y^r_{k})$, the hypergraphs $X_{k}|_{\left\{x^1_k, \ldots, x^{k - s + 2}_k\right\}}$ and $Y_{k}|_{\left\{y^1_k, \ldots, y^{k - s + 2}_k\right\}}$ are isomorphic. 
Let us show that for any $u_1, \ldots, u_l \in \{x^{k - s + 3}_k, \ldots, x^k_k\} \setminus V(\tilde X^1_{k} \cup \ldots \cup \tilde X^r_{k})$, where $l \ge 1$, there are no vertices $v_1, \ldots, v_{s - l} \in V(\tilde X^1_{k} \cup \ldots \cup \tilde X^r_{k})$
such that $\{u_1, \ldots, u_l, v_1, \ldots, v_{s - l}\} \in E(X_k)$. Suppose the contrary. We have $s - l \ge 2$. Let $v_1 \in \tilde X^{j_1}_k$, $v_2 \in \tilde X^{j_2}_k$. If $j_1 = j_2$, then the edge $\{u_1, \ldots, u_l, v_1, \ldots, v_{s - l}\}$ gives $1$-extension of $\tilde X^{j_1}_k$ which contradicts (II). If $j_1 \neq j_2$, then the existence of the edge $\{u_1, \ldots, u_l, v_1, \ldots, v_{s - l}\}$ implies that $d_{X_k}(\tilde X^{j_1}_k, \tilde X^{j_2}_k) \le 1$ which contradicts (III).
Therefore, $X_{k}|_{\{x^1_k, \ldots, x^k_k\}}$ has no edge containing at least one vertice from $V(X_k) \setminus V(\tilde X^1_{k} \cup \ldots \cup \tilde X^r_{k})$ (and the same is true for $Y_{k}|_{\{y^1_k, \ldots, y^k_k\}}$). 
Hence $X_{k}|_{\{x^1_k, \ldots, x^k_k\}}$ and $Y_{k}|_{\{y^1_k, \ldots, y^k_k\}}$ are isomorphic.\\

\textbf{Acknowledgements}
\\

This work has been supported by the RSF project 16-11-10014.

\renewcommand{\refname}{References}

\end{document}